\documentclass[hidelinks,onefignum,onetabnum]{siamart220329}
\usepackage[english]{babel}
\usepackage[margin=1.3in]{geometry}
\usepackage{braket,amsfonts}
\usepackage{amsmath, amssymb}

\usepackage{array}

\usepackage[caption=false]{subfig}
\newsiamthm{claim}{Claim}
\newsiamremark{remark}{Remark}
\newsiamremark{hypothesis}{Hypothesis}
\crefname{hypothesis}{Hypothesis}{Hypotheses}

\usepackage{algorithmic}
\Crefname{ALC@unique}{Line}{Lines}

\usepackage{graphicx,epstopdf}
\graphicspath{ {figs/} {./} }

\usepackage{amsopn}

\newcommand{\assign}{:=}

\newcommand{\mathd}{\mathrm{d}}

\newcommand{\tmop}[1]{\ensuremath{\operatorname{#1}}}

\newcommand{\tmtt}[1]{\texttt{#1}}

\newcounter{nntheorem}

\begin{document}

\title{Improved Laguerre Spectral Methods with Less Round-off Errors and Better Stability\thanks{\today,
\funding{This work was supported by the National Natural Science Foundation of China (Grant No. 12171467, 11771439)}}}

\author{
  Shenghe Huang\thanks{School of Mathematical Sciences, University of Chinese Academy of Sciences, Beijing  100049, China;
		NCMIS \& LSEC, Institute of Computational Mathematics and Scientific/Engineering Computing,
		Academy of Mathematics and Systems Science, Beijing 100190, China 
    		(\email{shenghe.huang@lsec.cc.ac.cn})
 	 }
  \and Haijun Yu\thanks{Corresponding author. 
		NCMIS \& LSEC, Institute of Computational Mathematics and Scientific/Engineering Computing,
		Academy of Mathematics and Systems Science, Beijing 100190, China;
		School of Mathematical Sciences, University of Chinese Academy of Sciences, Beijing  100049, China
		 (\email{hyu@lsec.cc.ac.cn})
	}
}

\maketitle

\begin{abstract}
  Laguerre polynomials are orthogonal polynomials defined on positive half
  line with respect to weight $e^{-x}$. They have wide applications in
  scientific and engineering computations. However, the exponential growth of Laguerre polynomials of high degree makes it hard to apply them to complicated systems that need to use large numbers of Laguerre bases. In this paper, we 
  introduce modified three-term recurrence formula to reduce the round-off error and to avoid overflow and underflow issues in generating generalized Laguerre polynomials and Laguerre functions. We apply the improved Laguerre methods to solve an elliptic equation defined on the half line. More than one thousand Laguerre bases are used in this application and meanwhile accuracy close to machine precision is achieved. The optimal scaling factor of Laguerre methods are studied and found to be independent of number of quadrature points in two cases that Laguerre methods have better convergence speeds than mapped Jacobi methods. 
\end{abstract}

% REQUIRED

% REQUIRED

\section{Introduction}

Spectral methods, which use orthogonal polynomials or orthogonal
trigonometric polynomials as bases, are a class of powerful methods to solve partial differential equations, especially when high accuracy is desired. For problems defined on half line $[0, +\infty)$, Laguerre polynomials are a natural choice, since they are orthogonal polynomials defined on half line with respect to weight $e^{-x}$. However, when apply Laguerre (resp. Hermite) polynomials to solve differential equations defined on half line (resp. whole line), there are several challenges~\cite{gottlieb_numerical_1977,funaro_computational_1990,shen_stable_2000}:
\begin{enumerate}
  \item Computations with Laguerre or Hermite polynomials usually lead to
  ill-conditioned algorithms and undesired round-off error instabilities.
  
  \item Direct using Laguerre or Hermite polynomials to approximate functions only provide meaningful results inside small intervals since the error estimation is obtained in weighted Sobolev spaces with fast decaying weights
  
  \item Standard Laguerre and Hermite polynomials have poor resolution properties when compared with other types of orthogonal polynomials.
\end{enumerate}
All these challenges have been studied and partially solved. For the first challenge, scaled Laguerre functions are introduced by Funaro\cite{funaro_computational_1990}. For the second challenge, it is proposed by Shen~\cite{shen_stable_2000} to use Laguerre functions as bases to form Laguerre spectral methods, this approach also (partially) avoids the first challenge. For the third challenge, it is found that the approximation by Hermite and Laguerre polynomials/functions can be significantly improved by choosing good scaling factors~\cite{tang_hermite_1993,shen_stable_2000}. Due to these pioneering works, the stability and efficiency of Laguerre and Hermite spectral methods are significant improved. They are widely used in scientific and engineering applications.
E.g. a stable Laguerre method is applied to solve nonlinear partial differential equations on a semi-infinite interval~\cite{guo_laguerregalerkin_2000}. Laguerre polynomials/functions are used to form composite spectral or spectral element methods for problems in unbounded domains~\cite{guo_composite_2001,xu_mixed_2002,shen_laguerre_2006a,wang_composite_2008,azaiez_laguerrelegendre_2009,chen_new_2012,shen_efficient_2016-1}. Laguerre methods have also been used to construct spectral numerical integrator for ordinary differential equations~\cite{guo_numerical_2007,guo_integration_2008}. Efficient time-splitting Laguerre-Hermite spectral methods are designed for Bose-Einstein condensates~\cite{bao_fourth-order_2005,bao_generalizedlaguerre_2008}. Laguerre and Hermite methods have also been used in optimal control problem~\cite{masoumnezhad_laguerre_2020} and fractional differential equations~\cite{bhrawy_efficient_2014,chen_laguerre_2018,tang_hermite_2018}. Scaling factors in different situations are further investigated~\cite{ma_hermite_2005,xia_efficient_2021}. A Laguerre method with scaling factor built in the Laguerre bases is designed and implemented with diagonalization technique~\cite{liu_spectral_2017,liu_fully_2017,zhang_efficient_2017}.

Despite these advances in Laguerre spectral methods for unbounded domains, the existing applications only used small numbers of Laguerre bases. This is partially due to the fact that standard subroutines to generate Laguerre polynomials and Laguerre functions suffer from overflow and underflow issues in generating large numbers of Laguerre-Gauss quadrature points. In this paper we introduce a simple but effective method to handle the underflow problem in calculating the values of Laguerre functions on large Gauss points. Numerical results show our method has better accuracy than existing methods given in ~\cite{funaro_computational_1990,shen_stable_2000}. On the same time, we find that the standard three-term recurrence formula to generate Laguerre polynomials and Laguerre functions leads to serious round-off errors for small arguments close to $0$. We introduce a modified 
three-term recurrence formula to reduce the round-off error. Numerical results show that the modified recurrence formula improve the accuracy of generating Laguerre polynomials by 2 to 4 decimal digits for small arguments close to 0, comparing to the standard recurrence formula. By using these improved basic algorithms for Laguerre spectral methods, we are able to solve a model elliptic equation defined on semi-infinite interval using more than one thousands Laguerre bases without convergence deterioration. To the best of our knowledge, these is the largest degree of Laguerre bases being used in Laguerre spectral methods. We also study the scaling factor of Laguerre spectral methods, it is found that for the two special cases that Laguerre methods convergent faster than mapped Jacobi methods, the optimal scaling factor is independent of number of Laguerre bases used. This means one can tune the scaling factor on coarse grids and use the optimal value on finer grids. We expect the improved basic algorithms for Laguerre methods and the findings on the optimal scaling factor will accelerate the development and applications of Laguerre spectral methods to large complex systems.

The remain part of this paper is organized as follows. In Section 2, we review some basic properties of generalized Laguerre polynomials/functions and related basic algorithms. In Section 3, we introduce the improved three-term recurrence formula with less round-off error for generalized Laguerre polynomials and functions, and the modified algorithm for generalized Laguerre functions to handle the underflow issue. An error estimation of the three-term recurrence formula is given. In Section 4, application of Laguerre spectral method to a model elliptic equation on semi-infinite interval is presented and the scaling factor is studied.  A short summary is given in Section 5.
\section{Preliminaries}

In this section, we review some basic properties of generalized Laguerre polynomials and functions, which are related to basic algorithms of Laguerre spectral methods. They can be found in standard reference books, e.g. \cite{abramowitz_handbook_1972,szego_orthogonal_1939, olver_nist_2010,canuto_spectral_2007,shen_spectral_2011}.
\subsection{Generalized Laguerre polynomials}

Generalized Laguerre polynomials, denoted by $\mathcal{L}_n^{(\alpha)}
(x)$ with ($\alpha > -1$) , are polynomials defined on the half line
$\mathbb{R}_+ \assign (0, +\infty)$ and orthogonal with respective to weight
\begin{equation}
  \omega_{\alpha} (x) = x^{\alpha} e^{- x}, \label{eq:wtg}
\end{equation}
i.e.
\begin{equation}
  \int_0^{+ \infty} \mathcal{L}_n^{(\alpha)} (x) \mathcal{L}_m^{(\alpha)} (x)
  \omega_{\alpha} (x) \mathd x = \gamma_n^{(\alpha)} \delta_{n m},
  \label{eq:orth}
\end{equation}
where $\delta_{n m}$ is the Kronecker delta symbol and
\begin{equation}
  \gamma_n^{(\alpha)} = \frac{\Gamma (n + \alpha + 1)}{n!} .
  \label{eq:normconst}
\end{equation}
In particular, $\mathcal{L}_n^{(0)} (x)$ is the usual Laguerre polynomials,
which will be denoted by $\mathcal{L}_n (x)$, its weight $\omega_0 (x)$ will
be denoted by $\omega (x) \assign e^{- x}$. We also have $\gamma_n^0 = 1$,
i.e. $\mathcal{L}_n (x)$ are normalized orthogonal bases in $L^2_{\omega}
(\mathbb{R}_+)$.

In addition to the orthogonality \eqref{eq:orth}, the following formulas
provide basic properties and algorithms of Laguerre polynomials. 
\begin{enumerate}
  \item Three-term recurrence formula
  \begin{equation}
    \left\{\begin{array}{l}
      \mathcal{L}^{(\alpha)}_0 (x) = 1, \qquad \mathcal{L}^{(\alpha)}_1 (x) = -
      x + \alpha + 1,\\
      (n + 1) \mathcal{L}_{n + 1}^{(\alpha)} (x) = (2 n + \alpha + 1 - x)
      \mathcal{L}_n^{(\alpha)} (x) - (n + \alpha) \mathcal{L}_{n -
      1}^{(\alpha)} (x), \quad n \geq 1.
    \end{array}\right. \label{eq:lapoly3t}
  \end{equation}
  From the recurrence formula, we obtain that
  \begin{equation}
    \mathcal{L}^{(\alpha)}_n (0) = \frac{\Gamma (n + \alpha + 1)}{n! \Gamma
    (\alpha + 1)} = \frac{\gamma_n^{(\alpha)}}{\Gamma (\alpha + 1)} .
    \label{eq:gla0}
  \end{equation}
  In particular, $\mathcal{L}_n (0) = 1$, $\mathcal{L}_n^{(1)} (0) = n + 1$.

  \item Some formulas related to the derivative of generalized Laguerre polynomials
  \begin{equation}
    \partial_x \mathcal{L}^{(\alpha)}_n (x) = - \mathcal{L}^{(\alpha + 1)}_{n
    - 1} (x) = - \sum_{k = 0}^{n - 1} \mathcal{L}^{(\alpha)}_k (x),
  \end{equation}
  \begin{equation}
    \partial_x \mathcal{L}^{(\alpha)}_{n + 1} (x) = \partial_x
    \mathcal{L}^{(\alpha)}_n (x) - \mathcal{L}^{(\alpha)}_n (x) .
    \label{eq:derivative}
  \end{equation}
  \item The scaling of Laguerre polynomials is given by Theorem 8.22.5 of \cite{szego_orthogonal_1939}.
  \begin{equation}
    \mathcal{L}^{(\alpha)}_n (x) = \frac{1}{\sqrt{\pi}} e^{\frac{x}{2}} x^{- \frac{\alpha}{2} - \frac14} n^{\frac{\alpha}{2} - \frac14} \times \Bigl[ \cos \bigl( 2 \sqrt{nx}
    - \frac{(2 \alpha + 1) \pi}{4} \bigr) + \frac{O (1)}{\sqrt{nx}}
    \Bigr], \label{eq:expgrowth}
  \end{equation}
  for $x \in [c / n, b]$ with $c$ and $b$ being fixed positive numbers. The
  exponential growth of Laguerre polynomials cause numerical  issue in standard (IEEE 754) floating-point arithmetic systems.
\end{enumerate}

\subsection{Generalized Laguerre functions}

Due to the exponential growth of the generalized Laguerre polynomials, the values $\mathcal{L}^{(\alpha)}_N (x)$ at Gauss quadrature points will exceed the representation limit of
IEEE 754 standard for double precision floating-point data with $N$ larger than several hundreds, thus cause overflow issue. To remove
this limitation, generalized \emph{Laguerre  functions} are considered, which
are defined as
\begin{equation}
  \widehat{\mathcal{L}}^{(\alpha)}_n (x) \assign e^{-x/2}
  \mathcal{L}^{(\alpha)}_n (x), \quad x \in \mathbb{R}_+, \quad \alpha > - 1.
  \label{eq:glafun}
\end{equation}
They are orthogonal with respective to weight
\begin{equation}
  \hat{\omega}_{\alpha} (x) = x^{\alpha}. \label{eq:glafwt}
\end{equation}
The corresponding formulas to generate generalized Laguerre functions and their derivatives are given below.
\begin{enumerate}
  \item Three-term recurrence formula:
  \begin{equation}
    \left\{\displaystyle 
      \begin{array}{l}
        \widehat{\mathcal{L}}^{(\alpha)}_0 (x) = e^{-x/2}, \qquad
        \widehat{\mathcal{L}}^{(\alpha)}_1 (x) = (\alpha + 1 - x) e^{-x/2}, \\   
    (n + 1) \widehat{\mathcal{L}}^{(\alpha)}_{n + 1} (x) = (2 n + \alpha + 1 -
    x) \widehat{\mathcal{L}}^{(\alpha)}_n (x) - (n + \alpha)
    \widehat{\mathcal{L}}^{(\alpha)}_{n - 1} (x), \quad n \geq 1.
      \end{array}
    \right.\label{eq:lafun3t}
  \end{equation}
  \item The recurrence formula for derivatives:
  \begin{equation}
    \begin{cases}
      \displaystyle \partial_x \widehat{\mathcal{L}}^{(\alpha)}_0 (x) = -\frac{1}{2} e^{- x
      / 2}, 
      \qquad \partial_x \widehat{\mathcal{L}}^{(\alpha)}_1 (x) = -
      \frac{\alpha +3-x}{2} e^{-x /2}, & \\ 
      \displaystyle \partial_x \widehat{\mathcal{L}}^{(\alpha)}_{n + 1} (x) = \partial_x
      \widehat{\mathcal{L}}^{(\alpha)}_n (x) - \frac{1}{2}
      \widehat{\mathcal{L}}^{(\alpha)}_n (x) - \frac{1}{2}
      \widehat{\mathcal{L}}^{(\alpha)}_{n + 1} (x), & n\ge 1.
     \end{cases} 
     \label{eq:lafunder}
  \end{equation}
\end{enumerate}

\subsection{Laguerre-Gauss type quadrature}

Let $\bigl\{ x_j^{(\alpha)},\, \omega_j^{(\alpha)} \bigr\}_{j = 0}^N$ be the set of
Laguerre-Gauss or Laguerre-Gauss-Radau quadrature nodes and weights.
\begin{itemize}
  \item For the Laguerre-Gauss quadrature,
   \begin{equation}
    \left\{ 
      \begin{array}{l}
        \bigl\{ x_j^{(\alpha)} \bigr\}_{j = 0}^N \text{ are the zeros of } \mathcal{L}^{(\alpha)}_{N + 1} (x), \text{ and} \\  
      \displaystyle \omega_j^{(\alpha)}  = \frac{\Gamma (N + \alpha + 1)}{(N + \alpha + 1)
        (N + 1) !} \frac{x_j^{(\alpha)}}{\bigl[ \mathcal{L}^{(\alpha)}_N
        (x_j^{(\alpha)}) \bigr]^2}, \quad 0 \leq j \leq N. 
      \end{array}
    \right.  \label{eq:lagswt}
     \end{equation}
  \item For the Laguerre-Gauss-Radau quadrature
  \begin{equation}
    \left\{ 
      \begin{array}{l}
      \displaystyle x_0^{(\alpha)} = 0, \quad  
      \bigl\{ x_j^{(\alpha)} \bigr\}_{j=1}^N \text{ are the zeros of } \partial_x 
      \mathcal{L}^{(\alpha)}_{N + 1} (x), \text{ and}\\ [2\jot]
      \displaystyle \omega_0^{(\alpha)}  = \frac{(\alpha + 1) \Gamma^2 (\alpha + 1)
       N!}{\Gamma (N + \alpha + 2)}, \\ [3\jot]
      \displaystyle \omega_j^{(\alpha)}  = \frac{\Gamma (N + \alpha + 1)}{N! (N + \alpha +
       1)} \frac{1}{[\mathcal{L}_N^{(\alpha)} (x_j^{(\alpha)})]^2}, \quad 1
       \leq j \leq N.
    \end{array}
    \right. \label{eq:lagsrwt}
  \end{equation}
\end{itemize}
With above nodes and weights, we have
\[ \int_0^{\infty} p (x) x^{\alpha} e^{- x} \mathd x = \sum_{j = 0}^N p
   (x_j^{(\alpha)}) \omega_j^{(\alpha)}, \quad \forall \, p \in P_{2 N +
   \delta}, \]
where $\delta = 1, 0$ for the Laguerre-Gauss quadrature and the
Laguerre-Gauss-Radau quadrature, respectively.

To form the Laguerre-Gauss and Laguerre-Gauss-Radau quadrature scheme, we need
find the zeros of the Laguerre polynomials, which can be obtained by calculating all 
eigenvalues of a banded matrix (see e.g. {\cite{golub_calculation_1969a}}{\cite{gautschi_orthogonal_1996}}). More precisely, by the
three-term recurrence formula \eqref{eq:lapoly3t}, the zeros $\{
x_j^{(\alpha)} \}_{j = 0}^N$ of $\mathcal{L}^{(\alpha)}_{N + 1} (x)$ are the
eigenvalues of the symmetric tridiagonal matrix (see e.g.
{\cite{shen_spectral_2011}}):
\begin{equation}
  A_{N + 1} = \left[\begin{array}{ccccc}
    a_0 & - \sqrt{b_1} &  &  & \\
    - \sqrt{b_1} & a_1 & - \sqrt{b_2} &  & \\
    & \vdots & \vdots & \vdots & \\
    &  & - \sqrt{b_{N - 1}} & a_{N - 1} & - \sqrt{b_N}\\
    &  &  & - \sqrt{b_N} & a_N
  \end{array}\right], \label{eq:3tEigenMatrix}
\end{equation}
where $a_j = 2 j + \alpha + 1$, $0 \leq j \leq N$, $b_j = j (j +
\alpha)$, $1 \leq j \leq N$. The approach of using this property to
calculate the zeros of orthogonal polynomials is referred to as the eigenvalue method.

For small values of $x_j^{(\alpha)}$, we have following asymptotic formula
\begin{equation}\label{eq:xGaussmin} 
  \sqrt{x_j^{(\alpha)}} = \frac{(j + 1) \pi + O (1)}{2 \sqrt{N + 1}}, \quad
   \text{for}\ x_j^{(\alpha)} \in (0, \eta), \ \alpha > - 1, 
  \end{equation}
where $\eta > 0$ is a fixed constant. For the largest point
\begin{equation} \label{eq:xN}
   x_N^{(\alpha)} = 4 N + 2 \alpha + 6 + O (N^{1 / 3}) . 
\end{equation}
By (6.31.13) of Szeg{\"o} (1975), we have also for $\alpha = 0$
\begin{equation}
  x_j = C_{j, N} \frac{(j + 2)^2}{N + 2}, \quad 0 \leq j \leq N,
  \quad \frac{1}{4} < C_{j, N} < 4.
\end{equation}
By the scaling \eqref{eq:xN}, \eqref{eq:expgrowth}, and the formula \eqref{eq:lagswt} and \eqref{eq:lagsrwt}, one easily see that the Laguerre-Gauss and Laguerre-Gauss-Radau quadrature weights decay exponentially for large $N$. E.g. the last weights for the 5
-point and 10-point Laguerre-Gauss formula are about $2.337\times 10^{-5}$, $9.912\times 10^{-13}$, respectively~\cite{olver_nist_2010}, and it is $6.770\times 10^{-23}$ for 16-point formula~\cite{shen_spectral_2011}. The unbalanced quadrature weights is not a good property for numerical computations. Fortunately, the Laguerre-Gauss quadrature can be directly extended to the Laguerre function approach, for which we define
\begin{equation}
  \hat{\omega}_j^{(\alpha)} : = e^{x_j^{(\alpha)}} \omega_j^{(\alpha)}, \quad
  0 \leq j \leq N. \label{eq:qwtfun}
\end{equation}
Then
\begin{equation}
  \int_0^{\infty} p (x) q (x) x^{\alpha} \mathd x = \sum_{j = 0}^N p
  (x_j^{(\alpha)}) q (x_j^{(\alpha)}) \hat{\omega}_j^{(\alpha)}, \quad \forall
  \, p \in \widehat{P}_K,\ q \in \widehat{P}_L, \; K + L \leq 2 N + \delta,
\end{equation}
where
\begin{equation}
  \widehat{P}_N \assign \bigl\{\phi : \phi = e^{-x/2} \psi, \quad \forall\, \psi \in
  P_N \bigr\}, \label{eq:PhatSpace}
\end{equation}
and $\delta=1, 0$ for the modified Laguerre-Gauss rule and the modified
Laguerre-Gauss-Radau rule, respectively.

\section{Improved Laguerre algorithms with less round-off error and better stability}

To implement a basic spectral method for solving differential equations, we
need basic quadrature and spectral transform algorithms, in which we need
first calculate the quadrature points and weights. The quadrature points can
be generated by using the eigenvalue method. After accurate quadrature points are
obtained, the quadrature weights can be obtained by using formula
\eqref{eq:lagswt}, \eqref{eq:lagsrwt} and \eqref{eq:qwtfun}, in which one
needs to evaluate the generalized Laguerre polynomials/functions at Gauss
quadrature points.

\subsection{Generating Gauss quadrature points}

\begin{figure}[tbhp]
  \centering
  \includegraphics[width=0.98\textwidth]{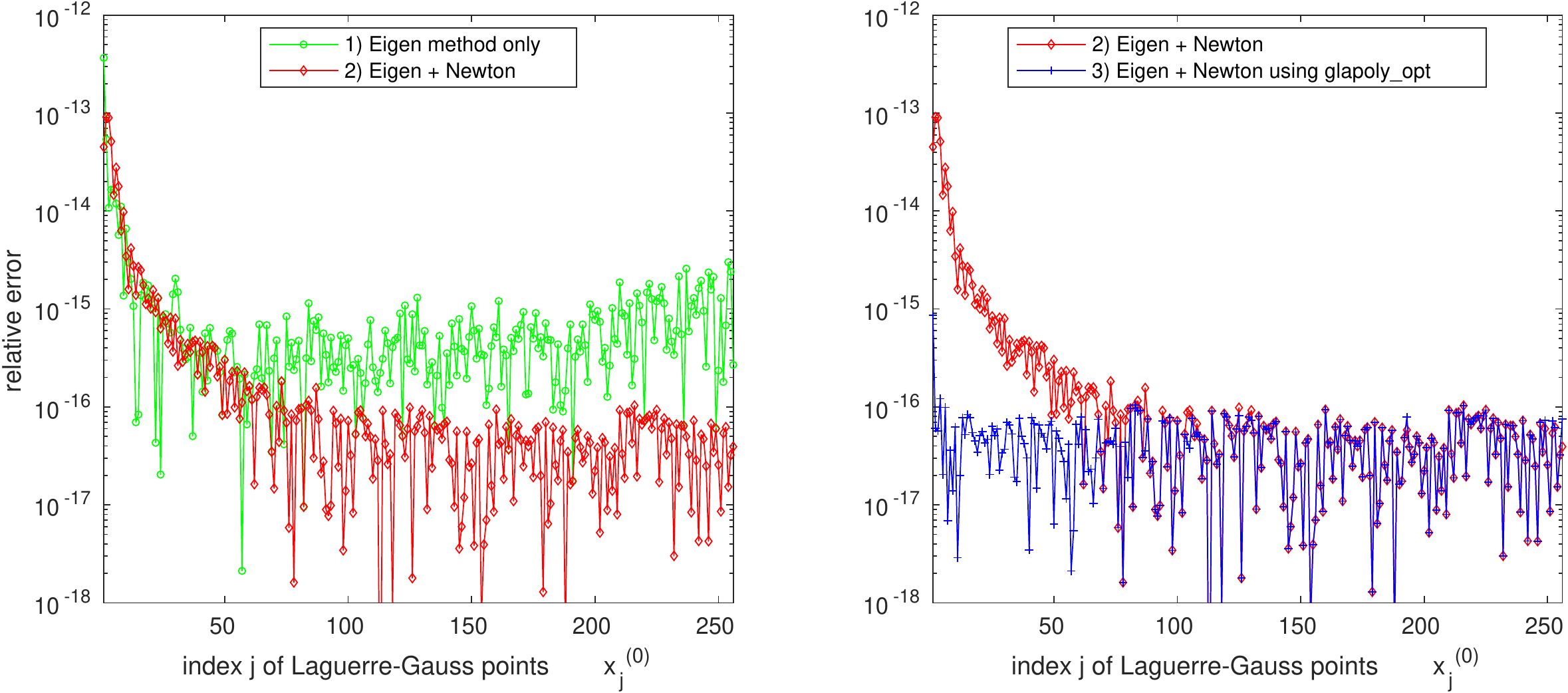}
  \caption{\label{fig:GaussPt} The relative errors of Laguerre-Gauss quadrature
  points~($N = 256$, $\alpha = 0$) calculated by different methods. In the
  right plot, {\tt{glapoly\_opt}} stands for generalized Laguerre polynomials
  generated using new approach \eqref{eq:glapoly3tnew}.}
\end{figure}

We first use eigenvalue method to generate the Laguerre-Gauss quadrature
points. Numerical results show that the eigenvalue method  can't
obtain full accuracy using floating-point arithmetic. Usually several significant
digits are lost~(see the left plot in Figure \ref{fig:GaussPt}). We then use the
following Newton's iterations to improve the accuracy of quadrature points
\begin{equation}
  x_j^{n + 1} = x_j^n - \frac{\mathcal{L}^{(\alpha)}_N (x_j^n)}{\partial_x
  \mathcal{L}^{(\alpha)}_N (x_j^n)}, \quad n = 0, 1, \ldots .
  \label{eq:NewtonM}
\end{equation}
To this end, we need a stable subroutine to calculate $\mathcal{L}^{(\alpha)}_N (x_j)$
and $\partial_x \mathcal{L}^{(\alpha)}_N (x_j)$ with high accuracy. To
evaluate the accuracy of the eigenvalue method and the classical three-term
recurrence formula to generate $\mathcal{L}^{(\alpha)}_N (x)$ and $\partial_x
\mathcal{L}^{(\alpha)}_N (x)$, we need a method to generate exact (reference)
solution, for which we use multiple precision arithmetic~(e.g.
{\tt{mpmath}} in Python or variable precision arithmetic (vpa) in Matlab
and Octave) to generate $\mathcal{L}^{(\alpha)}_N (x)$ and $\partial_x
\mathcal{L}^{(\alpha)}_N (x)$, and use them to refine the Gauss quadrature
points via \eqref{eq:NewtonM} to generate exact solutions. The results are
given in the left plot of Figure \ref{fig:GaussPt}, from which we see a clear
improvement of using Newton's iterations, especially for large $x_j$'s, where
the accuracy are improved more than 10 times and double precision limit is
reached. However, the improvements for small $x_j$'s are not
significant. The results for generating Laguerre-Gauss-Radau points are
similar, to save space they are not shown in the figure.

Theoretically, if $\mathcal{L}^{(\alpha)}_N (x)$ and $\partial_x
\mathcal{L}^{(\alpha)}_N (x)$ are evaluated accurately, the refinement of
Gauss points using Newton's iterations should be able to reach machine
accuracy. We think the big numerical errors for $x_j$'s close to 0 is mainly
caused by the round-off error of the coefficient $2 n + \alpha + 1 - x$ in the
three-term recurrence formula \eqref{eq:lapoly3t}. Since $x$ is small, but $2
n + \alpha + 1$ is not small, so some significant digits in $x$ will be lost
if we use floating-point system to store $2n +\alpha +1 -x$.

\subsection{Modified three-term recurrence algorithm}

To improve the accuracy of the three-term recurrence algorithm for small values of $x$, we take a new approach, by introducing
\[ \delta \mathcal{L}^{(\alpha)}_n = \mathcal{L}^{(\alpha)}_n -
   \mathcal{L}^{(\alpha)}_{n - 1}, \quad n \geq 1. \]
The new recurrence formula is given by
\begin{equation}
  \label{eq:glapoly3tnew} \left\{\begin{array}{ll}
    \displaystyle \mathcal{L}^{(\alpha)}_0 (x) = 1, 
    \qquad \mathcal{L}^{(\alpha)}_1 (x) = 1 + \alpha - x, 
    \qquad & \delta \mathcal{L}^{(\alpha)}_1 = \alpha - x, \\ [1\jot]
    \displaystyle \delta \mathcal{L}^{(\alpha)}_{n + 1} (x) = \frac{n + \alpha}{n + 1}
    \delta \mathcal{L}^{(\alpha)}_n (x) - \frac{x}{n + 1}
    \mathcal{L}^{(\alpha)}_n (x), & n \geq 1,\\ [2\jot]
    \displaystyle \mathcal{L}^{(\alpha)}_{n + 1} (x) = \mathcal{L}^{(\alpha)}_n (x) + \delta
    \mathcal{L}^{(\alpha)}_{n + 1} (x), & n \geq 1.
  \end{array}\right.
\end{equation}
If required, the derivatives $\partial_x \mathcal{L}^{(\alpha)}_{n + 1} (x)$ is calculated by using \eqref{eq:derivative} when $\mathcal{L}^{(\alpha)}_{n+1}(x)$ is obtained.

Numerical results of generated Laguerre-Gauss points $x_j^{(\alpha)}$ for
$\alpha = 0$, $N = 256$ by using the standard and new three-term recurrence
formulas are shown in the right plot of Figure \ref{fig:GaussPt}, from which
we see that the new approach perform significantly better than the old standard
approach. The only point with relative
error larger than machine accuracy is the first quadrature point
next to 0. But it still has an improvement of 2 significant digits over the standard recurrence formula with Newton's iterations.

\begin{figure}[htbp]
  \centering
  \includegraphics[width=0.47\textwidth]{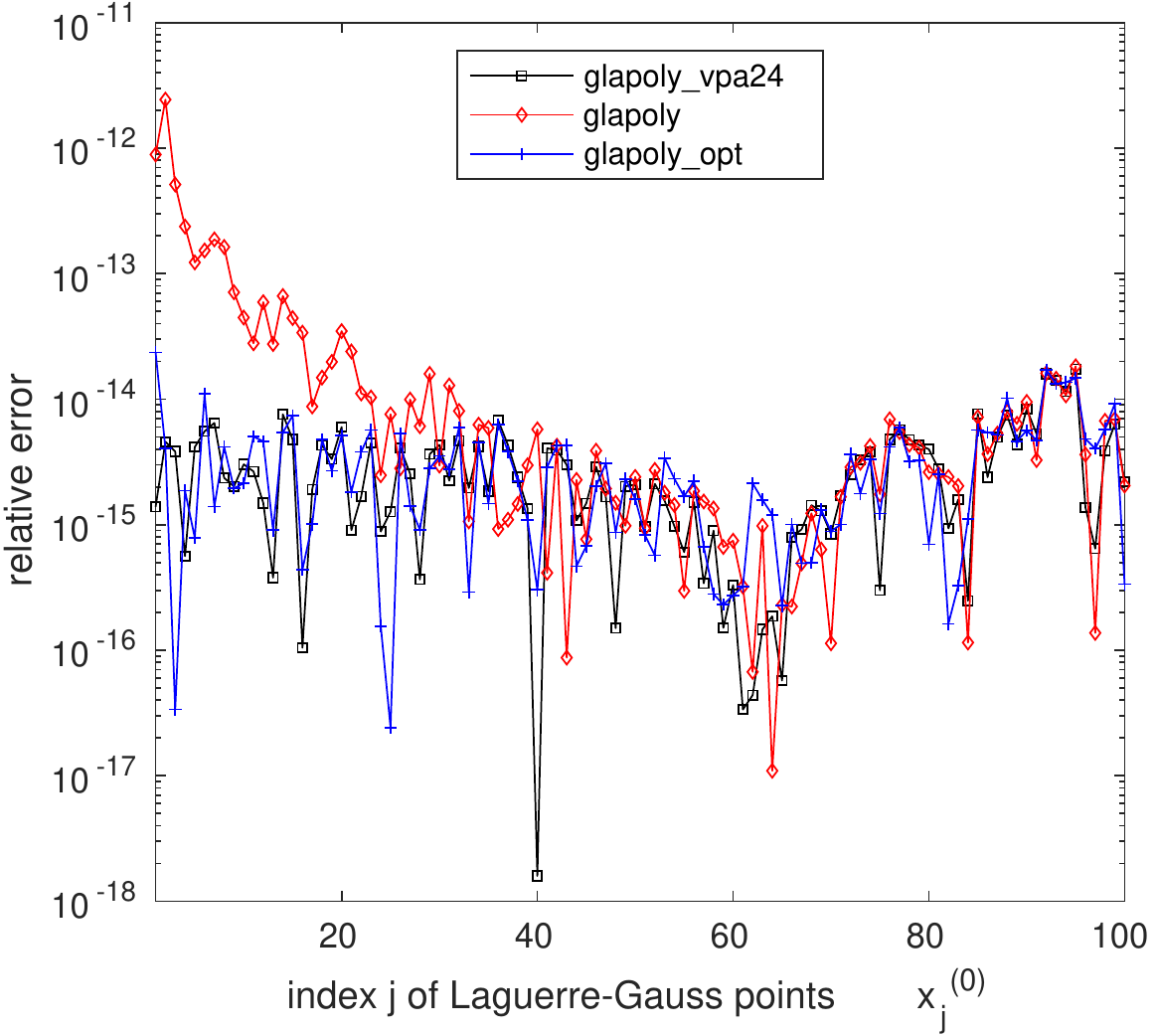}
  \includegraphics[width=0.47\textwidth]{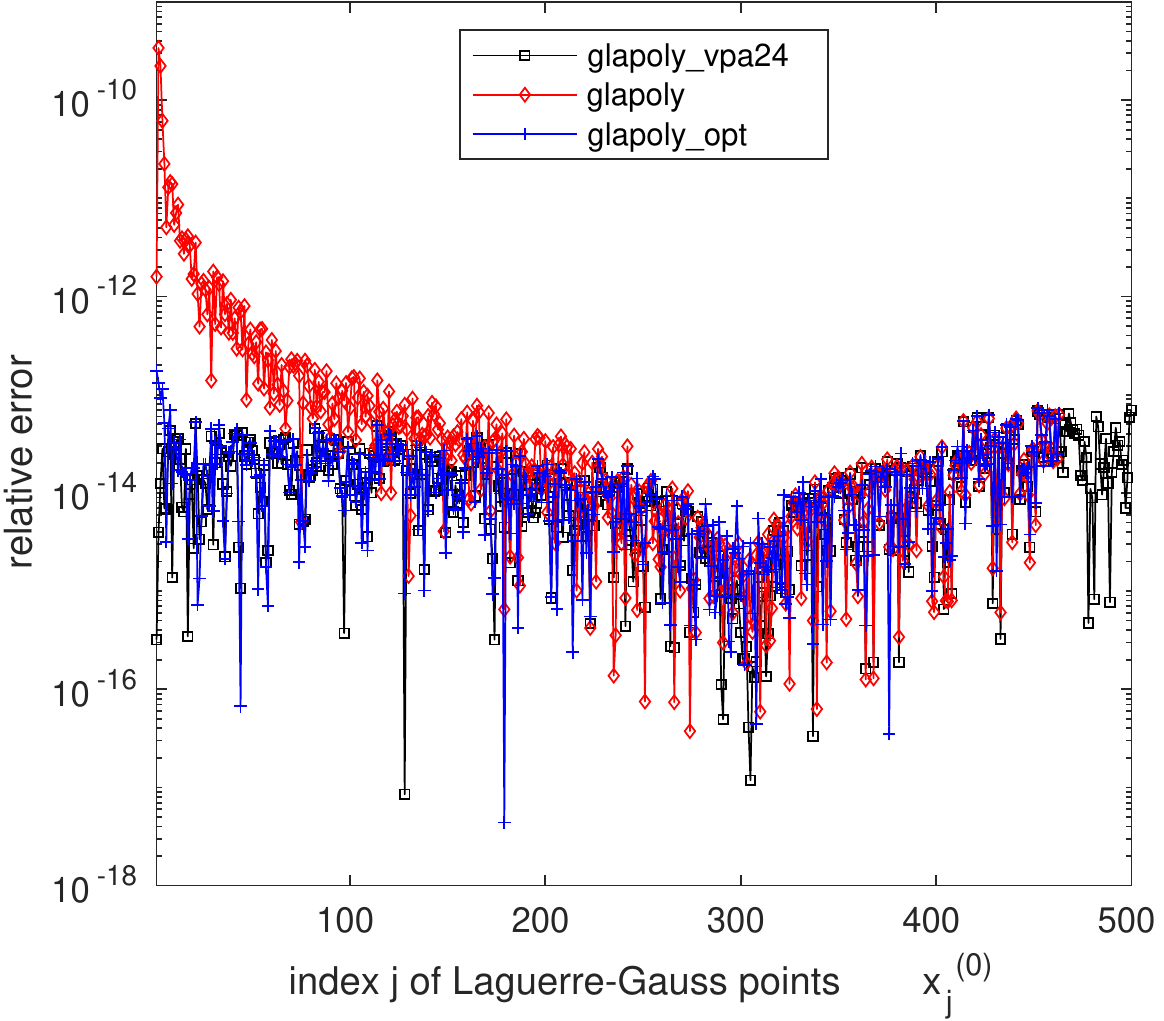}
  \caption{\label{fig:lapoly} $\mathcal{L}^{(\alpha)}_{N - 1}
  (x_j^{(\alpha)})$ for $\alpha = 0$ calculated by two methods: glapoly()
  stands for the traditional formula \eqref{eq:lapoly3t}, {\tt glapoly\_opt}
  stands for the new approach \eqref{eq:glapoly3tnew}, {\tt glapoly\_vpa24} is a
  reference solution obtained by using variable-precision algorithms of 24
  digits. Left: $N = 100$; Right: $N = 500$. For glapoly and glapoly\_opt in
  the right plot, the last 37 values are {\tt{NaN}} and not plotted.}
\end{figure}

We then test the error of classical procedures to evaluate Laguerre
polynomials at given Gauss points, which are used in the formulas of
Gauss quadrature weights \eqref{eq:lagswt} and \eqref{eq:lagsrwt}. The results
of using two methods to generate $\mathcal{L}^{(\alpha)}_{N - 1}
(x_j^{(\alpha)}), j = 1, \ldots, N$ for $\alpha = 0, N = 100$ are presented in
the left plot of Figure \ref{fig:lapoly}, where {\tt{glapoly}} stands
for the standard three-term recurrence algorithm, {\tt{glapoly\_vpa24}}
stands for results generated by standard three-term recurrence algorithm using
24-digit floating-point arithmetic, {\tt{glapoly\_opt}} stands for the results
using new recurrence formula \eqref{eq:glapoly3tnew}. We observe that the
classical three-term recurrence algorithm is not very accurate, several
significant digits are lost especially for small $x_j$'s. On the other hand,
{\tt{glapoly\_opt}} gives much accurate results, especially for smaller
$x_j$'s.

We next test the results of generating Laguerre polynomials for a larger $N$.
In the right plot of Figure \ref{fig:lapoly}, we present the results of
generating $\mathcal{L}^{(\alpha)}_{N - 1}$ at $N$ Laguerre-Gauss points using
different methods with $N = 500$. We find that with data stored in
double precision, the Laguerre polynomial $\mathcal{L}^{(\alpha)}_{499} (x_j^{(\alpha)})$
exceeds the representing region of double-precision for $j>463$. We tested
that using double-precision floating-point number system to obtain \eqref{eq:lagswt} with $a = 0, 1, 2, 3, 4, 5, 6$, the maximum Laguerre-Gauss quadrature scheme can be generated is $N =
362$. We will propose a method to overcome the overflow
problem later in this section, but in next subsection, we first give an error
estimate of the three-term recurrence formulas.

\subsection{An error estimate for the three-term recurrence algorithms}

For three-term recurrence formula \eqref{eq:lapoly3t}, we have
\[ \mathcal{L}_{n + 1}^{(\alpha)} (1 + \delta_{n + 1}) = \frac{2 n + \alpha +
   1 - x (1 + \eta_x)}{n + 1} (1 + \eta^a_n) \mathcal{L}_n^{(\alpha)} (1 +
   \delta_n) - \frac{n + \alpha}{n + 1} (1 + \eta^b_n) \mathcal{L}_{n -
   1}^{(\alpha)} (1 + \delta_{n - 1}), \]
where $\eta_x, \eta^a_n, \eta_n^b$ account for the relative round-off errors
of the corresponding coefficients, $\delta_{n + 1}, \delta_n$, $\delta_{n -
1}$ are the relative errors of $\mathcal{L}^{(\alpha)}_{n + 1}$,
$\mathcal{L}^{(\alpha)}_n$, $\mathcal{L}^{(\alpha)}_{n - 1}$, respectively.
Using notation $e_n = \mathcal{L}^{(\alpha)}_n (x) \delta_n$, we have
\begin{eqnarray*}
  e_{n + 1} & = & \frac{2 n + \alpha + 1 - x (1 + \eta_x)}{n + 1} (1 +
  \eta^a_n) (\mathcal{L}_n^{(\alpha)} + e_n) - \frac{n + \alpha}{n + 1} (1 +
  \eta^b_n) (\mathcal{L}_{n - 1}^{(\alpha)} + e_{n - 1})\\
  &  & - \frac{2 n + \alpha + 1 - x}{n + 1} \mathcal{L}^{(\alpha)}_n +
  \frac{n + \alpha}{n + 1} \mathcal{L}^{(\alpha)}_{n - 1}\\
  & = & -\frac{x}{n + 1} \eta_x (1 + \eta^a_n) (\mathcal{L}_n^{(\alpha)} +
  e_n) + \frac{2 n + \alpha + 1 - x}{n + 1} \eta^a_n (\mathcal{L}_n^{(\alpha)}
  + e_n) + \frac{2 n + \alpha + 1 - x}{n + 1} e_n\\
  &  & - \frac{n + \alpha}{n + 1} \eta^b_n (\mathcal{L}_{n - 1}^{(\alpha)} +
  e_{n - 1}) - \frac{n + \alpha}{n + 1} e_{n - 1}\\
  & = & -\frac{x}{n + 1} \mathcal{L}_n^{(\alpha)} \eta_x + \frac{2 n + \alpha
  + 1 - x}{n + 1} \mathcal{L}_n^{(\alpha)} \eta^a_n + \frac{2 n + \alpha + 1 -
  x}{n + 1} e_n\\
  &  & - \frac{n + \alpha}{n + 1} \mathcal{L}_{n - 1}^{(\alpha)} \eta^b_n -
  \frac{n + \alpha}{n + 1} e_{n - 1} + \xi_n,
\end{eqnarray*}
where $\xi_n$ includes high-order error terms, which can be bounded by
$\varepsilon (4 | e_n | + | e_{n-1} |)$.  Here $\varepsilon$ is the bound of
relative round-off error of the floating-point system used. For double precision data in IEEE 754 standard, $\varepsilon \approx 2.22\times 10^{-16}$. The above result
can be written as
\begin{equation}
  \left\{\begin{array}{l}
    e_{n + 1} = \frac{2 n + \alpha + 1 - x}{n + 1} e_n - \frac{n + \alpha}{n +
    1} e_{n - 1} + \zeta_n, \quad n \geq 1,\\
    e_0 = 0,\\
    e_1 = \mathcal{L}_1^{(\alpha)} \eta_0^a, \quad |\eta_0^a| \sim \varepsilon .
  \end{array}\right. \label{eq:err3t}
\end{equation}
where
\begin{equation}
  \zeta_n = \xi_n - \frac{x}{n + 1} \mathcal{L}_n^{(\alpha)} \eta_x + \frac{2
  n + \alpha + 1 - x}{n + 1} \mathcal{L}_n^{(\alpha)} \eta^a_n - \frac{n +
  \alpha}{n + 1} \mathcal{L}_{n - 1}^{(\alpha)} \eta^b_n . \label{eq:zn1}
\end{equation}
For $x \ll 1$, $n\gg 1$, $\zeta_n \sim \Bigl( 2 + \dfrac{x}{n + 1} \Bigr) |
\mathcal{L}_n^{(\alpha)} | \varepsilon + | \mathcal{L}_{n - 1}^{(\alpha)} |
\varepsilon$.

\begin{theorem*}
  Assume that $\alpha > -1$, Then, for any $\eta > 0$ with $3 -\alpha
  -x -\eta > 0$ , the error propagation given by \eqref{eq:err3t} satisfies
  the following estimate:
  \[ E_{n+1} 
     \leq \left\{\begin{array}{ll}
       E_1 + \frac{(n+1)(n+2)(2n+3)}{6 \eta}  \max_s \zeta_s^2, 
       & \text{if}\ 1 -2\alpha -x -\eta \geq 0,\\
       \beta_n E_1 + \frac{(n-3)(n+1)^2+29 \beta_n}{\eta} \max_s \zeta_s^2, 
       & \text{if}\ -1.5<1 - 2\alpha -x -\eta < 0\ \text{and}\ \alpha\ge 0,
     \end{array}\right. \]
  for $n \geq 1$, where $E_n \assign xe_{n}^2 + (n+ \alpha) (e_{n} - e_{n-1})^2$ and $\beta_n = \frac{\Gamma (n+2+\alpha) / \Gamma
  (2+\alpha)}{\Gamma (n+3 -\alpha -x -\eta) / \Gamma (3 -\alpha -x -\eta)}$.

  \begin{proof}
    We first rewrite the first equation of \eqref{eq:err3t} as
    \[ (e_{n + 1} - e_n) - \frac{n + \alpha}{n + 1} (e_n - e_{n - 1}) = -
       \frac{x}{n + 1} e_n + \zeta_n . \]
    Then multiply $d_{n + 1} : = e_{n + 1} - e_n$ on both sides, we obtain
    \[ \frac{1 - \alpha}{n + 1} d_{n + 1}^2 + \frac{1}{2} \frac{n + \alpha}{n
       + 1} [d_{n + 1}^2 - d_n^2 + (d_{n + 1} - d_n)^2] = \frac{1}{2}
       \frac{x}{n + 1} [e_n^2 - e_{n + 1}^2 + d_{n + 1}^2] + \zeta_n d_{n +
       1}, \]
    i.e.
    \[ \frac{x}{n + 1} [e_{n + 1}^2 - e_n^2] + \frac{n + 2 - \alpha - x}{n +
       1} d_{n + 1}^2 - \frac{n + \alpha}{n + 1} d_n^2 + \frac{n + \alpha}{n +
       1} (d_{n + 1} - d_n)^2 = 2 \zeta_n d_{n + 1}, \]
    or
    \[ xe_{n + 1}^2 + (n + 2 - \alpha - x) d_{n + 1}^2 + (n + \alpha) (d_{n +
       1} - d_n)^2 = xe_n^2 + (n + \alpha) d_n^2 + 2 (n + 1) \zeta_n d_{n + 1}
       . \]
    Then, by using Cauchy inequality with $\eta$ on $\zeta_n d_{n + 1}$: $2
    \zeta_n d_{n + 1} \leq \frac{n + 1}{\eta} \zeta_n^2 + \frac{\eta}{n +
    1} d_{n + 1}^2$, we get
    \begin{equation}
      x e_{n+1}^2 + (n +2 -\alpha -x -\eta) d_{n+1}^2 
      + (n+\alpha)(d_{n+1} - d_n)^2 \leq x e_n^2 + (n+\alpha) d_n^2 
      + \frac{(n+1)^2}{\eta} \zeta_n^2.
    \end{equation}
    Define
    \[ r_n \assign \max \left\{ 1, \frac{n + 1 + \alpha}{n + 2 - \alpha - x -
       \eta} \right\}. \]
    Then for $n \geq 1$, $3 - \alpha - x - y > 0$,
    \[ E_{n + 1} \leq r_n E_n + r_n \frac{(n + 1)^2}{\eta} \zeta_n^2 . \]
    Then by mathematical induction, we obtain
    \begin{equation} \label{eq:basic-estimate}
      E_{n + 1} \leq (\Pi_{k = 1}^n r_k) E_1 + \sum_{s = 1}^n (\Pi_{k =
      s}^n r_k) \frac{(s + 1)^2}{\eta} \zeta_s^2, \quad n \geq 1.
    \end{equation}
    When $1 - 2\alpha -x - \eta \geq 0$, we have $r_n = 1$, thus
    \begin{equation} \label{eq:est1}
      E_{n+1} \leq E_1 + \sum_{s=1}^n \frac{(s+1)^2}{\eta} \zeta_s^2 
    \leq E_1 + \frac{(n+1)(n+2)(2n+3)}{6 \eta} \max_s \zeta_s^2.
    \end{equation}
    {When $-1.5<1 - 2\alpha -x -\eta < 0$, $r_n = \frac{n+1 +\alpha}{n+2-\alpha -x -\eta}>1$. Let $b_s = (\Pi_{k=s}^n r_k) \frac{(s+1)^2}{\eta}$, then
    \[ b_s / b_{s-1} = \frac{(s+1)^2}{s^2} \frac{s+1-\alpha-x-\eta}{s+\alpha} 
    = \left( 1+ \frac{2}{s} + \frac{1}{s^2}\right) \left(1-\frac{\beta}{s+\alpha}\right), \]
    where $0<\beta = 2\alpha + x + \eta - 1 \leq 1.5$, which means 
    \begin{equation}
      \frac{b_s}{b_{s-1}} = 1 + \frac{(2-\beta)s + 2\alpha - 2\beta}{s(s+\alpha)} + \frac{1}{s^2} \frac{s+\alpha-\beta}{s+\alpha}\ge 1\quad
      \text{if}\ s\ge 5, \alpha \ge 0.
    \end{equation}
    Define $\beta_n^s = \Pi_{k=s}^n r_k$, and $\beta_n = \beta_n^1$. From above inequality and \eqref{eq:basic-estimate}, we have 
    \begin{align} 
      E_{n + 1} & \leq \beta_n E_1 + 
    \left( \frac{(n-3)(n+1)^2}{\eta}+ \sum_{s=1}^3 \beta_n^s \frac{(s+1)^2}{\eta} \right) \max_s \zeta_s^2 \nonumber \\
    & \le \beta_n E_1 + 
    \left( \frac{(n-3)(n+1)^2}{\eta}+ \beta_n  \frac{29}{\eta} \right) \max_s \zeta_s^2.  \label{eq:est2}
    \end{align}
    Combine the two cases \eqref{eq:est1} and \eqref{eq:est2}, we obtain the desired estimate.}
  \end{proof}
\end{theorem*}

From above theorem, it is easy to deduce the following result.

\begin{theorem}
  For $- 1 < \alpha \leq 1 / 4, x < 1 / 4$, taking $\eta = 1 / 4$, we
  have error estimate for \eqref{eq:err3t}:
  \begin{align}
    | e_{n + 1} | 
    &\leq \sqrt{\frac{1}{| x |} \biggl( (1 + \alpha + x) e_1^2 +
    \frac{2(n+1)(n+2)(2n+3)}{3} \max_s \zeta_s^2 \biggr) } \\
    & \leq \left( 2 + \frac{2(n+2)^{3/2}}{\sqrt{3}} \right) \frac{\max \{ | e_1 |, | \zeta_s | \}}{| x |^{1 / 2}} . \label{eq:err3tc1}
  \end{align}
  For $\alpha \ge 0$, $3 -\alpha -x -\eta > 0$, 
  $-1.5\le 1 - 2\alpha -x -\eta < 0$,
  we have
  \begin{eqnarray*}
    | e_{n + 1} | & \leq & \sqrt{\frac{1}{| x |} \left[ \beta_n (1 +
    \alpha + x) e_1^2 + \frac{29\beta_n + (n-3)(n+1)^2}{\eta} \max_s \zeta_s^2 \right]}\\
    & \leq & \left( 2 \sqrt{\beta_n} + \frac{5.5\sqrt{\beta_n}+n^{1 / 2} (n + 1)
    }{\sqrt{\eta}} \right) \frac{\max \{ | e_1 |, | \zeta_s | \}}{| x |^{1 /
    2}},
  \end{eqnarray*}
  Here $\beta_n \sim O ((n +1 +\alpha)^{2 \alpha +x  +\eta -1})$. In
  particular, for $\alpha = 1$, $x < 1 / 4$, taking $\eta = 1 / 4$, we have
  \[ | e_{n + 1} | \leq (O (n^{3 / 4}) + O (n^{3 / 2}) ) \frac{\max \{ |
     e_1 |, | \zeta_s | \}}{| x |^{1 / 2}} \sim O (n^{3 / 2}) \frac{\max \{ |
     e_1 |, | \zeta_s | \}}{| x |^{1 / 2}} . \]
\end{theorem}

\begin{remark}
  Since $e_1 \sim | \mathcal{L}_1^{(\alpha)} | \varepsilon$, $\zeta_s \sim
  \left( \left( 2 + \frac{x}{n + 1} \right) | \mathcal{L}_s^{(\alpha)} | + |
  \mathcal{L}_{s - 1}^{(\alpha)} | \right) \varepsilon$. For $x$ close to the
  smallest Gauss point, by \eqref{eq:xGaussmin}, we have $x \sim O(1/N)$, thus  $| e_{N} | \lesssim O (N^2 |
  \mathcal{L}_n^{(\alpha)} | \varepsilon)$. Note that this is a worst-case error estimate. In practice, the numerical errors will be smaller due to the cancellations in summation of random round-off errors. 
\end{remark}

\begin{figure}[htbp]
  \centering
  \includegraphics[width=0.47\textwidth]{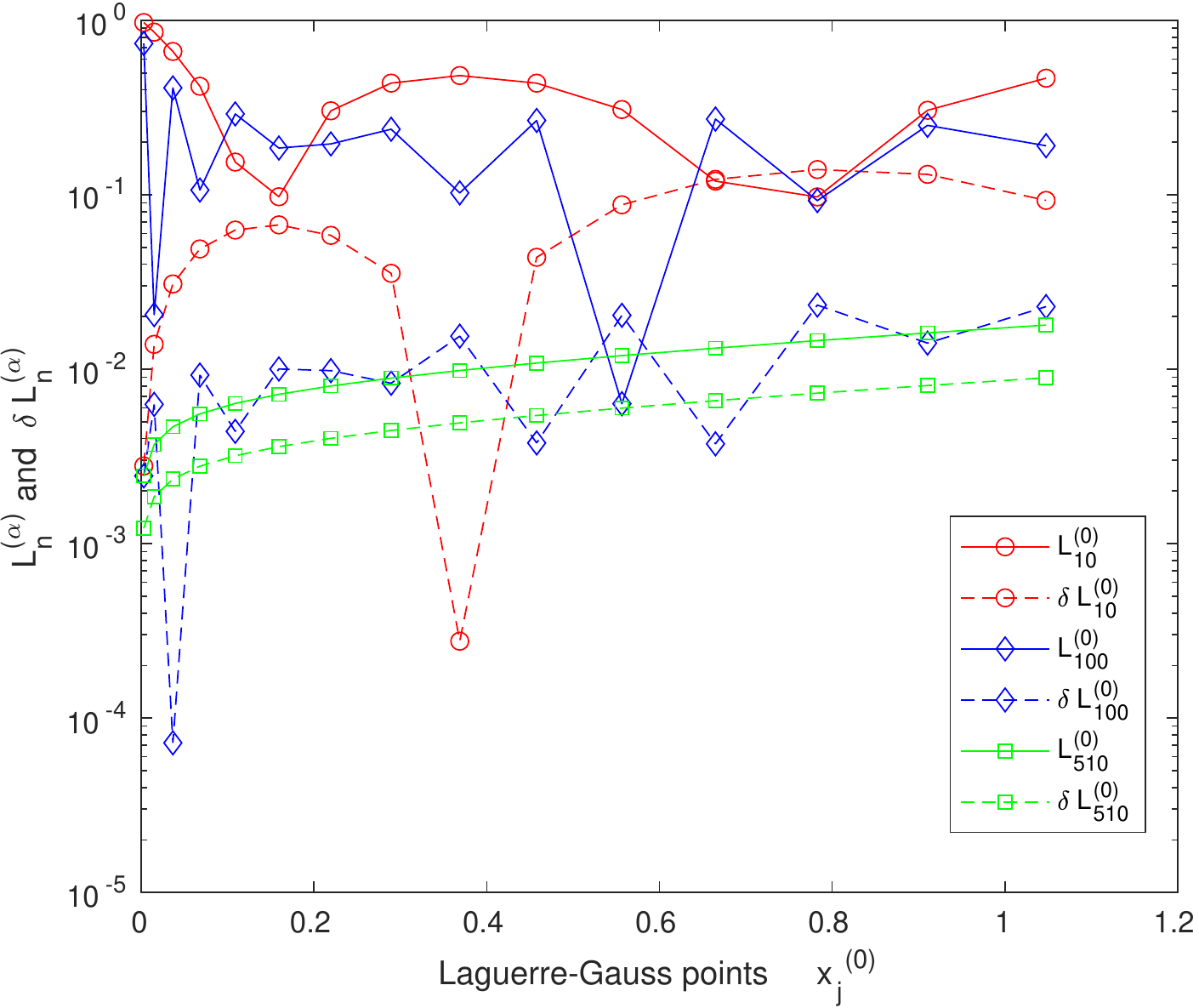}
  \includegraphics[width=0.47\textwidth]{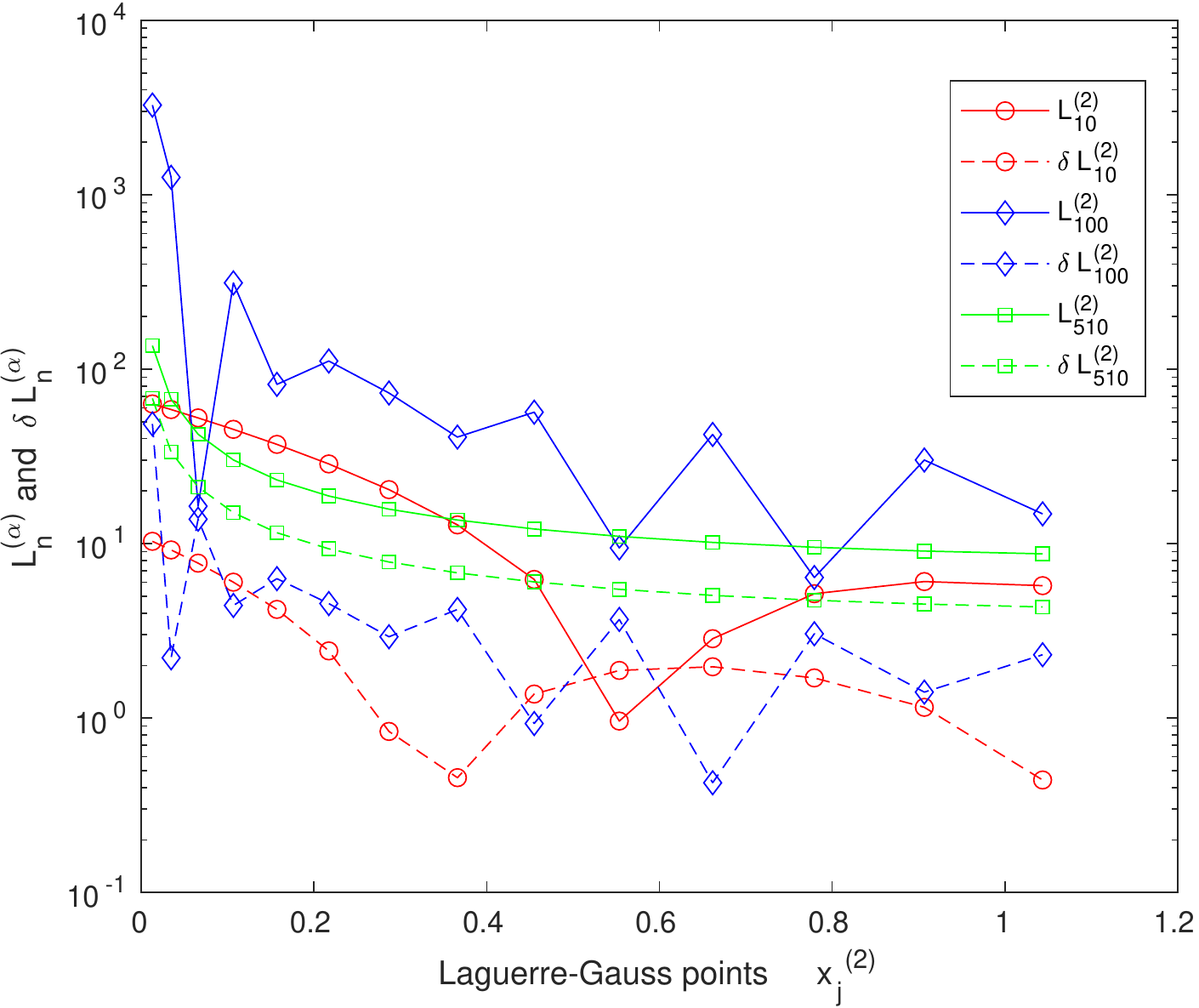}
  \caption{\label{fig:LdL} $|\mathcal{L}^{(\alpha)}_n (x_j^{(\alpha)})|$  and
  $| \delta \mathcal{L}^{(\alpha)}_n (x_j^{(\alpha)}) |$, $n = 10,
  100, 510$, $\alpha = 0$ (left) and $\alpha = 2$ (right). $x_j^{(\alpha)}$
  are from 512th order Laguerre-Gauss points.}
\end{figure}

For iteration \eqref{eq:glapoly3tnew}, we have
\[ \left\{\begin{array}{l}
     \delta \mathcal{L}^{(\alpha)}_{n + 1} (x) (1 + \delta_{n + 1}^a) =
     \frac{n + \alpha}{n + 1} (1 + \eta_n^a) \delta \mathcal{L}^{(\alpha)}_n
     (1 + \delta_n^a) - \frac{x}{n + 1} (1 + \eta_n^b)
     \mathcal{L}^{(\alpha)}_n (x) (1 + \delta_n^b),\\ [2\jot]
     \mathcal{L}^{(\alpha)}_{n + 1} (1 + \delta_{n + 1}^b) =
     \mathcal{L}^{(\alpha)}_n (x) (1 + \delta_n^b) + \delta
     \mathcal{L}^{(\alpha)}_{n + 1} (x) (1 + \delta_{n + 1}^a) .
   \end{array}\right. \]
From above equation, we obtain error propagation equation
\begin{equation} 
  \left\{\begin{array}{l}
     d_{n + 1} = \frac{n + \alpha}{n + 1} d_n - \frac{x}{n + 1} e_n +
     \zeta_n^{\delta}, \quad n \geq 1,\\
     e_{n + 1} = e_n + d_{n + 1}, \quad n \geq 1,\\
     e_0 = 0, \quad d_1 = e_1 .
   \end{array}\right.   
   \label{eq:new3terrp}
\end{equation}
where
\begin{equation}
  \zeta_n^{\delta} = \frac{n + \alpha}{n + 1} \delta \mathcal{L}^{(\alpha)}_n
  \eta_n^a - \frac{x}{n + 1} \mathcal{L}^{(\alpha)}_n \eta_n^b +
  \xi^{\delta}_n, 
  \qquad \xi^{\delta}_n = \frac{n + \alpha}{n + 1} \eta_n^a d_n
  - \frac{x}{n + 1} \eta_n^b e_n . \label{eq:zn2}
\end{equation}
Substituting for $d_{n+1} = e_{n+1} - e_n$, and $d_n = e_n - e_{n-1}$,
the first equation in \eqref{eq:new3terrp} becomes
\begin{equation}
  e_{n + 1} - e_n = \frac{n + \alpha}{n + 1} (e_n - e_{n - 1}) - \frac{x}{n +
  1} e_n + \zeta_n^{\delta} . \label{eq:err3tnew}
\end{equation}
We notice that the error equation is equivalent to \eqref{eq:err3t} except the
difference in $\zeta_n$ and $\zeta_n^{\delta}$. From equation \eqref{eq:zn2},
we know than $\zeta_n^{\delta} \sim O \left( \bigl| \delta
\mathcal{L}^{(\alpha)}_n \bigr| + \frac{x}{n + 1} \bigl|
\mathcal{L}^{(\alpha)}_n \bigr| \right) \varepsilon$. We observe that $\bigl| \delta
\mathcal{L}^{(\alpha)}_n (x) \bigr|$ is in general smaller than
$\bigl|\mathcal{L}^{(\alpha)}_n (x)\bigr|$  for $x$ small (see Figure \ref{fig:LdL}),
thus $\zeta_n^{\delta}$ is in general smaller than $\zeta_n$ for $x \ll n+1$. So it is expected that the iteration algorithm \eqref{eq:glapoly3tnew}
produce smaller numerical errors than \eqref{eq:lapoly3t} for $x$ not too
large, which has been demonstrated in Figure \ref{fig:lapoly}.

\subsection{Stable algorithms for Laguerre functions and Gauss points}

To remove the overflow problem associated with Laguerre polynomials approach,
we show how to use the Laguerre function approach to generate corresponding
Gauss quadrature scheme.

\begin{figure}[htbp]
  \centering
  \includegraphics[width=0.49\textwidth]{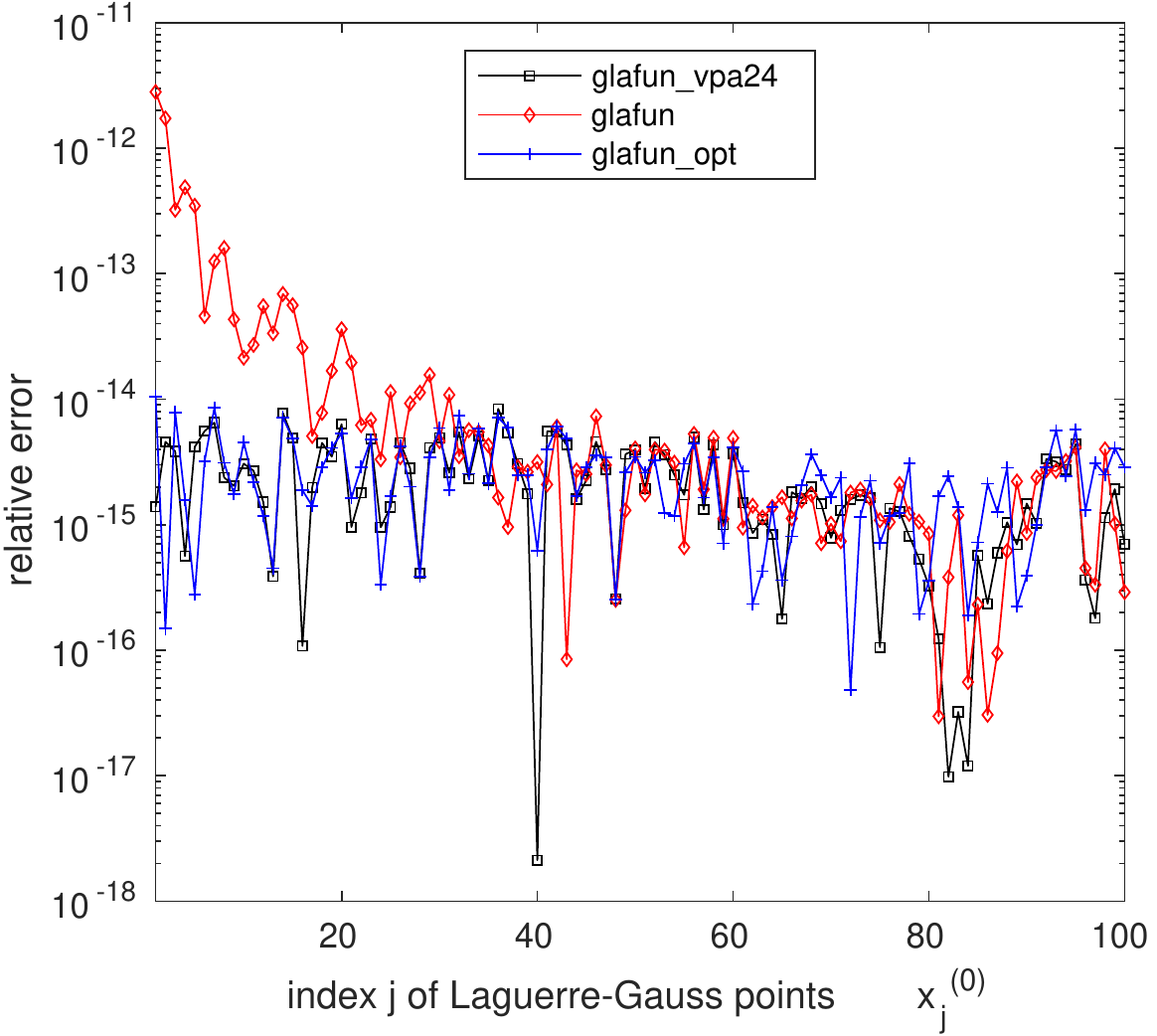}
  \includegraphics[width=0.49\textwidth]{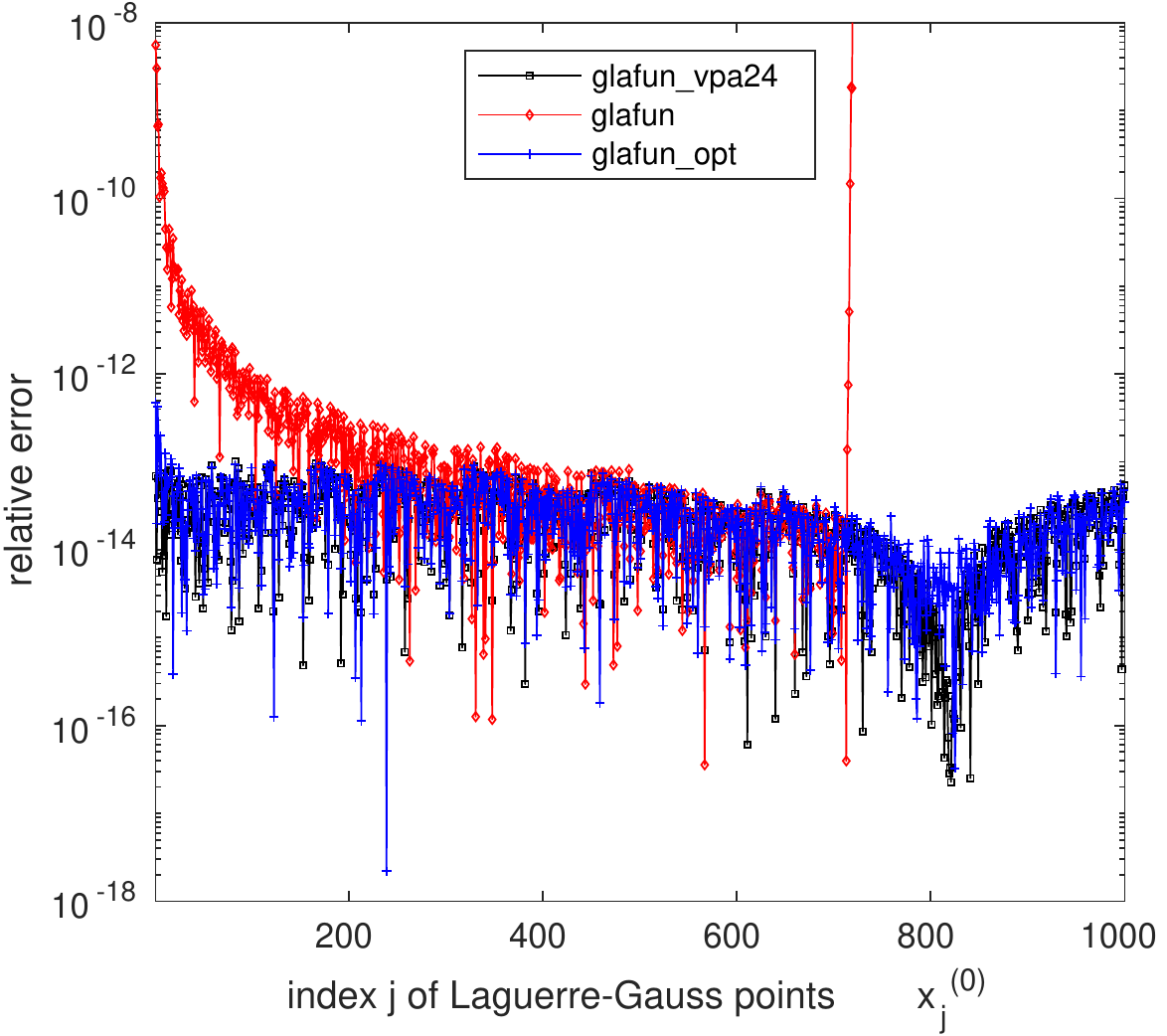}
  \caption{\label{fig:glafun}The relative error of
  $\widehat{\mathcal{L}}^{(\alpha)}_{N - 1} (x_j^{(\alpha)})$ calculated using
  different methods. Left: $\alpha = 0, N = 100$, Right: $\alpha = 0, N =
  1000$.}
\end{figure}

We first use a procedure similar to \eqref{eq:glapoly3tnew} to generate
generalized Laguerre functions.
\begin{equation}
  \label{eq:glaf3term} \left\{\begin{array}{ll}
    \widehat{\mathcal{L}}^{(\alpha)}_0 (x) = e^{-x/2}, 
    \qquad
    \widehat{\mathcal{L}}^{(\alpha)}_1 (x) = (1 + \alpha - x) e^{-x/2},
    \qquad &\delta \widehat{\mathcal{L}}^{(\alpha)}_1 =(\alpha - x) e^{-x/2},\\ [1\jot]
    \displaystyle \delta \widehat{\mathcal{L}}^{(\alpha)}_{n + 1} (x) = \frac{n + \alpha}{n
    + 1} \delta \widehat{\mathcal{L}}^{(\alpha)}_n (x) - \frac{x}{n + 1}
    \widehat{\mathcal{L}}^{(\alpha)}_n (x), & n \geq 1,\\ [3\jot]
    \widehat{\mathcal{L}}^{(\alpha)}_{n + 1} (x) =
    \widehat{\mathcal{L}}^{(\alpha)}_n (x) + \delta
    \widehat{\mathcal{L}}^{(\alpha)}_{n + 1} (x), & n \geq 1. 
  \end{array}\right.
\end{equation}
The results of calculating $\widehat{\mathcal{L}}^{(\alpha)}_{N - 1}
(x_j^{(\alpha)})$ for $\alpha = 0, N = 200$ using several methods are
presented in the left plot of Figure \ref{fig:glafun}. We see the new
recurrence formula \eqref{eq:glaf3term}, which is denoted by
{\tt{glafun\_opt}}, improve the results of standard approach (denoted by
{\tt{glafun}}) significantly, especially for small quadrature points. It is
almost as accurate as the results obtained using 24-digit high-precision
arithmetic (denoted by {\tt{glafun\_vpa24}}), which means the major
numerical errors are from the floating-point round-off error of the Gauss
points.

However, to generate the Laguerre-Gauss quadrature rule, if we
directly multiply the whole weight $e^{- x_j^{(\alpha)}/2}$ to the Laguerre
polynomials to obtain the corresponding Laguerre functions, we will encounter
double precision underflow for $N \geq 383$ when $\alpha = 0, 1$, since
for $x \gtrsim 745$, $e^{- x}$ will be treated as zero in double
precision system of IEEE 754 standard.

There are several methods have been proposed to overcome the overflow and underflow problem,
e.g. the scaled Laguerre function method proposed by Funaro
{\cite{funaro_computational_1990}}, and the sign and amplitude factor method
suggested by Shen {\cite{shen_stable_2000}}. Here, we propose to multiply the
weights gradually and adaptively. More precisely, in the iteration procedure,
we monitor the magnitude of the iteration variable, if it is larger than a
given threshold, then we multiply it by a portion of $e^{- x_j^{(\alpha)} /
2}$ to make it smaller than another threshold. At the end of the iteration, we
apply the remaining parts of the weight to obtain desired results. The overall
procedure is given in \cref{alg:glafun-opt}.

\begin{algorithm}
  \caption{\tt{glafun\_opt}}
  \label{alg:glafun-opt}
  \begin{algorithmic}[1]
  \STATE{\textbf{Input:} $N$, $\alpha$, $x$, $K_1 = 32$, $K_2 = 32$}
  \STATE{\textbf{Output:} $y = \widehat{\mathcal{L}}^{(\alpha)}_N (x)$ }

  \IF{$N = 0$}
	\STATE set $y=\exp(-x/2)$ 
	\RETURN $y$ 
  \ENDIF 
  \IF{$N = 1$}
	\STATE set $y = (1+\alpha - x) \exp(-x / 2)$
	\RETURN $y$ 
  \ENDIF 
  \STATE{set $L = 1 + \alpha - x$;}
  \STATE{set $\mathd L = \alpha - x$;}
  \STATE{set $x_b = x / 2$;}
  \FOR{k=1 to $N - 1$}
    \STATE set $\mathd L = ((k + \alpha) \mathd L - xL) / (k + 1)$;
    \STATE set $L = L + \mathd L$;

    \IF{$k = 1$ or $| L | > \exp (K_1)$}
    	\STATE set $x_c = \min (\max (\log | L | + K_2, 0), x_b)$;
    	\STATE set $\mathd L = \mathd L \exp (- x_c)$;
    	\STATE set $L = L \exp (- x_c)$;
    	\STATE set $x_b = x_b - x_c$;
    \ENDIF
  \ENDFOR

    %\STATE set $\mathd y = \mathd L \exp (- x_b)$;
    \STATE set $y = L \exp (- x_b)$;

  \RETURN $y$ %, $\mathd y$
  \end{algorithmic}
\end{algorithm}

Note that in the \cref{alg:glafun-opt}, $K_1$, $K_2$ are two
given positive constants. We require $K_1 + K_2 < 80$ to avoid overflow when
$x$ and $L$ are stored in single or double precision. In all the numerical
experiments presented in this paper, we set $K_1 = K_2 = 32$. The results of
$N = 1000$ is given in the right plot of Figure \ref{fig:glafun}. \ We see the
algorithm {\tmtt{glafun\_opt}} is as accurate as the {\tmtt{glafun\_vpa24}}
for large $x_j$'s, while the standard 3-term recurrence formula is only
accurate for $500 < j \lesssim 700$, but overflow for large $x_j$'s. For small
$x_j$'s, the results of {\tmtt{glafun\_opt}} is slight worse than the
{\tmtt{glafun\_vpa24}} results, but still has improvement up to 4 significant
digits comparing to the standard {\tmtt{glafun}} method.

\begin{figure}[htbp]
  \centering
  \includegraphics[width=0.49\textwidth]{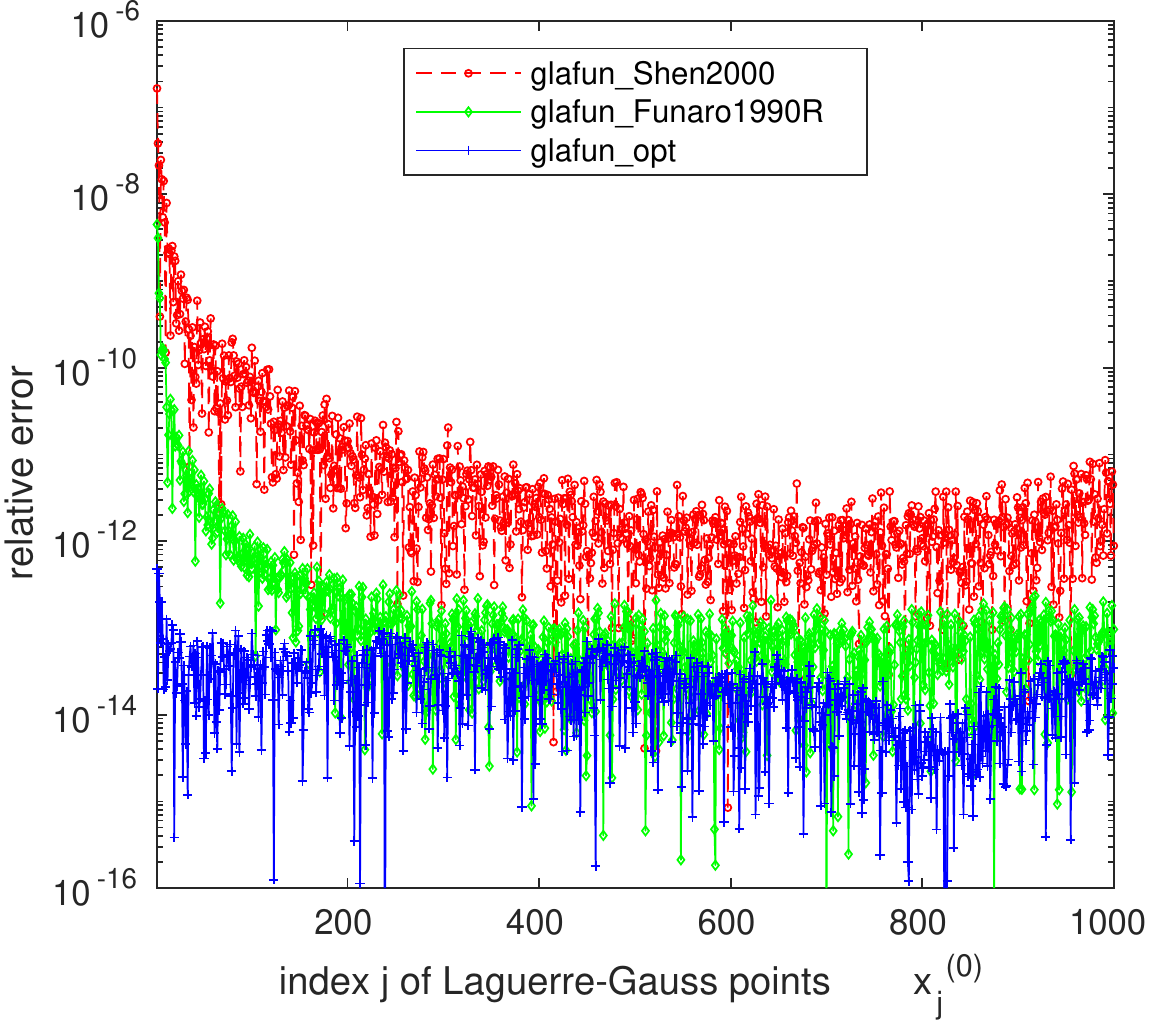}
  \includegraphics[width=0.49\textwidth]{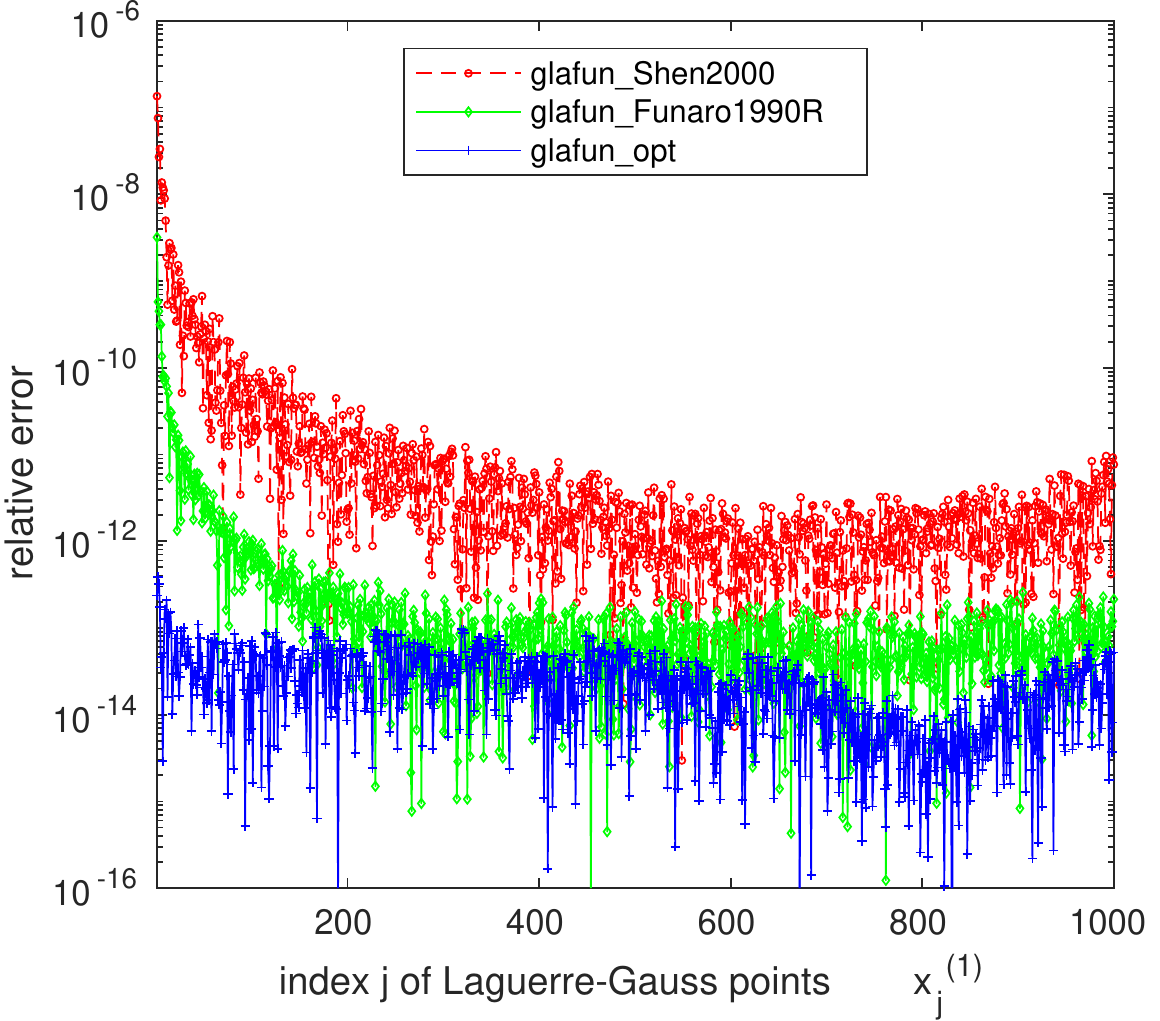}
  \caption{\label{fig:glafun_cmp}The relative error of
  $\widehat{\mathcal{L}}^{(\alpha)}_{N - 1} (x_j^{(\alpha)})$ calculated using
  different methods. Left: $\alpha = 0, N = 1000$, Right: $\alpha = 1, N =
  1000$.}
\end{figure}

Before ending this section, we compare the new recurrence algorithm to
generate generalized Laguerre functions with the method proposed in
{\cite{shen_stable_2000}} and the method proposed in
{\cite{funaro_computational_1990}}. Note that in
{\cite{funaro_computational_1990}}, {\emph{scaled}} Laguerre functions are
introduced, which is not equivalent to standard Laguerre functions. But we can
obtain standard Laguerre functions by multiplying an appropriate factor on the
scaled Laguerre functions. The numerical results are given in Figure
\ref{fig:glafun_cmp}, from which we see all the three methods can generate
generalized Laguerre functions for $N > 500$ with reasonable accuracy, but the
method proposed here give significantly better results than the other two
methods.

\section{Application: Laguerre method with best scaling factor}

\subsection{The Laguerre spectral method for an elliptic equation defined on half line}

Consider the following model equation
\begin{equation}
  - u_{x \, x} + \gamma u = f, \quad x \in \mathbb{R}_+, \quad \gamma > 0;
  \quad u (0) = 0, \quad \lim_{x \rightarrow +\infty} u(x) = 0.
  \label{eq:elleq}
\end{equation}
The corresponding weak formulation is: find $u \in H_0^1 (\mathbb{R}_+)$, such that
\begin{equation}
  a(u, v) \assign (u', v') + \gamma (u, v) = (f, v), \quad \forall\, v \in
  H_0^1 (\mathbb{R}_+), \label{eq:elleqwk}
\end{equation}
for $f \in H_0^1 (\mathbb{R}_+)'$. Here $H_0^1 (\mathbb{R}_+)'$ is the dual space of $H_0^1 (\mathbb{R}_+) \assign \bigl\{ v \in H^1(\mathbb{R}_+)\;\text{and}\; v(0)=0 \bigr\}$.

It is clear than the problem admits a
unique solution, since
\[ a (u, u) = | u |_1^2 + \gamma \| u \|^2 \geqslant \min (1, \gamma) \| u
   \|_1^2, \quad \forall \, u \in H_0^1 (\mathbb{R}_+) . \]
Define
\[ \widehat{P}_N^0 = \tmop{span} \Set{ \psi_n (x)
   \assign e^{-x/2} \bigl(\mathcal{L}^{(0)}_n (x) - \mathcal{L}^{(0)}_{n+1} (x) \bigr), \: n = 0, \ldots, N-1 } .
\]
The Laguerre spectral-Galerkin approximation to \eqref{eq:elleqwk} is : Find
$u_N \in \widehat{P}_N^0$, such that
\begin{equation}
  a(u_N, v_N) = (\hat{I}_M  f, v_N), \quad \forall \, v_N \in \widehat{P}_N^0 .
  \label{eq:LagGal}
\end{equation}
Here $\hat{I}_M$ denotes the interpolation from $C (\mathbb{R}_+)$ to
$\hat{P}_M$ using $M+1$ Laguerre-Gauss points.

The error of the above Laguerre-Galerkin method \eqref{eq:LagGal} is given by
{\cite{shen_spectral_2011}}:

\begin{theorem}
  Consider $\alpha = 0$, \ $\gamma > 0$. If $u \in H_0^1 (\mathbb{R}_+)$,
  $\hat{\partial}_x u \assign \left( \partial_x + \frac{1}{2} \right) u \in
  \hat{B}_0^{m - 1} (\mathbb{R}_+)$, $f \in C (\bar{\mathbb{R}}_+) \cap
  \hat{B}_0^k (\mathbb{R}_+)$ and $\hat{\partial}_x f \in \hat{B}_0^{k - 1}
  (\mathbb{R}_+)$ with $1 \leqslant m \leqslant N+1$, $1 \leq  k \leq M + 1$, then we have
  \begin{multline}
    \| u - u_N \|_1 \leqslant c \sqrt{\frac{(N - m + 1) !}{N!}} \|
    \hat{\partial}_x^m u \|_{\hat{\omega}_{m - 1}} \\
    + c \sqrt{\frac{(M - k + 1)
    !}{M!}} \Bigl(\| \hat{\partial}_x^k f \|_{\hat{\omega}_{k - 1}} + (\ln M)^{1 /
    2} \| \hat{\partial}_x^k f \|_{\hat{\omega}_k}\Bigr), \label{eq:elleqerr}
  \end{multline}
\end{theorem}

where $c$ is a positive constant independent of $m, k, N, M, u$ and $f$. Here,
the weighted Sobolev space $\hat{B}^m_{\alpha} (\mathbb{R}_+)$ together with
its associated norm are defined as
\[ \hat{B}_{\alpha}^m (\mathbb{R}_+) \assign \Set{ u \mid \hat{\partial}_x^k u
   \in L_{\hat{\omega} _{\alpha + k}}^2 (\mathbb{R}_+), \; 0 \leqslant k
   \leq m }, 
   \qquad 
   \| u \|_{\hat{B}_{\alpha}^m} \assign \biggl( \sum_{k=0}^m
   \| \hat{\partial}_x^k u \|_{\hat{\omega} _{\alpha + k}}^2 \biggr)^{1/2} .
\]
We will test the Laguerre method for \eqref{eq:elleq} with exponential and
algebraic decay solutions. In particular, we will test three cases with
exact solutions given below:
\begin{equation}
  u_1 (x) = \sin (kx) e^{-x}, 
  \qquad u_2 (x) = (1 + x)^{-r}, 
  \qquad u_3 (x) =  \cos (kx) (1 + x)^{-r} . \label{eq:utest}
\end{equation}

\subsection{The scaling factor}

It is known that the accuracy of spectral methods for unbounded domain can be
improved significantly by using a good scaling factor
{\cite{boyd_asymptotic_1984}}{\cite{tang_hermite_1993}}{\cite{shen_recent_2009}}.
An effective empirical rule is proposed for the Hermite spectral method by
Tang~\cite{tang_hermite_1993}, where a scaling factor $\beta$ is added to
the Hermite function basis $\psi_k (\beta x)$, such that the largest Hermite-Gauss quadrature point $\max_j | x^{(N)}_j |$ after the scaling is located at
a given threshold position
\begin{equation}
  \max_j | x_j^{(N)} | / \beta = x_M, \label{eq:tangbeta}
\end{equation}
where $x_M$ is determined by the criterion that the solution $u (x)$ is
numerically neglectable for $|x| > x_M$, i.e. for a given tolerance error
$\varepsilon$,
\[ | u (x) | < \varepsilon, \quad \forall\: | x | \geqslant x_M . \]
This approach is very effective and successful. Similar treatments have been
applied to Laguerre spectral methods
{\cite{shen_stable_2000}}{\cite{shen_recent_2009}}. However, by using this
approach, one has to determine $x_M$, and the choice $\beta$ given by
\eqref{eq:tangbeta} depends on the number of quadrature points used. We found
that for special exponential decay functions, e.g. the $u_1 (x)$ defined in
\eqref{eq:utest}, an optimal scaling factor $\beta$ can be determined by a
rigorous error estimation. To see this, we first apply the linear transform
$y = \beta x$ to equation \eqref{eq:elleq} to get
\begin{equation}
  - v_{y y} + \frac{\gamma}{\beta^2} v = \frac{1}{\beta ^2} g (y), \quad v (0)
  = 0, \quad \lim_{y \rightarrow + \infty} v (y) = 0, \label{eq:elleqmap}
\end{equation}
where $v (y) = u (x) = u (y / \beta)$, $g (y) = f (y / \beta)$. It is obvious
that for any $u^{\ast} (x)$ satisfies \eqref{eq:elleq}, $v^{\ast} (y) =
u^{\ast} (y / \beta)$ satisfies \eqref{eq:elleqmap}. The difference is that
the convergence speed using $\hat{P}_N^0$ to approximate problem
\eqref{eq:elleq} and \eqref{eq:elleqmap} are different. By applying
\eqref{eq:elleqerr} to \eqref{eq:elleqmap}, and omitting the interpolation
error for simplicity (suppose we choose $M > N$ to make the interpolation
error smaller than the projection error), we obtain
\begin{equation} 
  \| v - v_N \|_1 {\leqslant c_1}  \sqrt{\frac{(N - m + 1) !}{N!}} \|
   \hat{\partial}_y^m v \|_{\hat{\omega}_{m-1}} 
   \label{eq:elleqerrmap}
\end{equation}
where $c_1 = c \frac{\max (1, \gamma / \beta^2)}{\min (1, \gamma / \beta^2)}$.

Obviously, the convergence speed depends on the choice of measurement. We should
measure error using certain norm of $u (x)$. Notice that
\[ \| v - v_N \|^2 = \int_0^{\infty} | v (y) - v_N (y) |^{^2} \mathd y = \beta
   \int_0^{\infty} | u (x) - u_N (x) |^2 \mathd x = \beta \| u - u_N \|^2, \]
\[ | v - v_N |_1^2 = \int_0^{\infty} | \partial_y v (y) - \partial_y v_N (y)
   |^{^2} \mathd y = \frac{1}{\beta} \int_0^{\infty} | \partial_x u -
   \partial_x u_N (x) |^2 \mathd x = \frac{1}{\beta} | u - u_N |_1^2 . \]

We have
\begin{equation}
  | u - u_N |_1 = \sqrt{\beta} | v - v_N |_1 \leqslant c\sqrt{\beta} \frac{\max
  (1, \gamma / \beta^2)}{\min (1, \gamma / \beta^2)} \sqrt{\frac{(N - m + 1)
  !}{N!}} \| \hat{\partial}_y^m v \|_{\hat{\omega}_{m - 1}}
  \label{eq:errH1semi}
\end{equation}
\begin{equation}
  \| u - u_N \| = \frac{1}{\sqrt{\beta}} \| v - v_N \| \leq
  \frac{c}{\sqrt{\beta}} \frac{\max (1, \gamma / \beta^2)}{\min (1, \gamma /
  \beta^2)} \sqrt{\frac{(N - m + 1) !}{N!}} \| \hat{\partial}_y^m v
  \|_{\hat{\omega}_{m - 1}} . \label{eq:errL2}
\end{equation}
The errors depend on $\| \hat{\partial}_y^m v \|_{\hat{\omega}_{m - 1}}$,
which is also related to $\beta$, since
\begin{equation}
  \| \hat{\partial}_y^m v \|_{\hat{\omega}_{m - 1}}^2 = \int_0^{\infty} |
  (\partial_y + 1 / 2)^m u (y / \beta) |^2 y^{m - 1} \mathd y = \beta^m
  \int_0^{\infty} | (\partial_x / \beta + 1 / 2)^m u (x) |^2 x^{m - 1} \mathd
  x. \label{eq:beta-err}
\end{equation}
When $\beta \rightarrow 0$, we have
\[ \beta^{\pm \frac{1}{2}} \frac{\max (1, \gamma / \beta^2)}{\min (1, \gamma /
   \beta^2)} \rightarrow O \left( \frac{\gamma}{\beta^{2 \mp 1 / 2}} \right),
   \quad \tmop{and} \quad \| \hat{\partial}_y^m v \|_{\hat{\omega}_{m - 1}}^2
   \rightarrow \frac{1}{\beta^m} \int_0^{\infty} | \partial_x^m u (x) |^2 x^{m
   - 1} \mathd x. \]
When $\beta \rightarrow + \infty$, we have
\[ \beta^{\pm \frac{1}{2}} \frac{\max (1, \gamma / \beta^2)}{\min (1, \gamma /
   \beta^2)} \rightarrow O \left( \frac{\beta^{2 \pm 1 / 2}}{\gamma} \right),
   \quad \tmop{and} \quad \| \hat{\partial}_y^m v \|_{\hat{\omega}_{m - 1}}^2
   \rightarrow \beta^m \int_0^{\infty} \left| \frac{1}{2^m} u (x) \right|^2
   x^{m - 1} \mathd x. \]
So, for $m \geqslant 1$, $\int_0^{\infty} \left| \frac{1}{2^m} u (x) \right|^2
x^{m - 1} \mathd x > 0$ and $\int_0^{\infty} | \partial_x^m u (x) |^2 x^{m -
1} \mathd x > 0$, the right hand sides of \eqref{eq:errH1semi} and
\eqref{eq:errL2} have unbounded limit at $\beta \rightarrow 0$ and $\beta
\rightarrow \infty$. Therefore, there exists an optimal value $\beta$ leads to
optimal smallest value of the right hand sides of \eqref{eq:errH1semi} and
\eqref{eq:errL2}.

\subsubsection{Optimal scaling factor for the exponential decay function}

Now, we consider a special case $u (x) = e^{zx},$ where $z \in \mathbb{C}$,
$\Re (z) < 0$. Then
\begin{eqnarray*}
  \| \hat{\partial}_y^m v \|_{\hat{\omega}_{m - 1}}^2 & = & \beta^m
  \int_0^{\infty} | (\partial_x / \beta + 1 / 2)^m e^{zx} |^2 x^{m - 1} \mathd
  x\\
  & = & \beta^m \int_0^{\infty} \bigl| (z / \beta + 1 / 2)^{2 m} \bigr|  | e^{2 zx} |
  x^{m - 1} \mathd x\\
  & = & \biggl| \left( z / \sqrt{\beta} + \sqrt{\beta} / 2 \right)^{2 m}
  \biggr| \int_0^{\infty}  | e^{2 zx} |  x^{m - 1} \mathd x\\
  & = & \left| z / \sqrt{\beta} + \sqrt{\beta} / 2 \right|^{2 m} \frac{(m -
  1) !}{\Re (2z)^m}
\end{eqnarray*}
To make $\| \hat{\partial}_y^m v \|_{\hat{\omega}_{m - 1}}^2$ small, we need
make $\left| z / \sqrt{\beta} + \sqrt{\beta} / 2 \right|$ small. Suppose $z =
a + bi$, then
\[ \left| z / \sqrt{\beta} + \sqrt{\beta} / 2 \right| = \sqrt{\left(
   \frac{a}{\sqrt{\beta}} + \frac{\sqrt{\beta}}{2} \right)^2 +
   \frac{b^2}{\beta}} = \sqrt{\frac{a^2 + b^2}{\beta} + \frac{\beta}{4} + a}
\]
The minimal value is obtained when 
\begin{equation}  
  \beta = \beta_{\ast} = 2 \sqrt{a^2 + b^2} = 2 | z |, \label{eq:beta-exp}
\end{equation}
and
\begin{equation}
  \min_{\beta} \| \hat{\partial}_y^m v \|_{\hat{\omega}_{m - 1}}^2 = \bigl(| z |
  +\mathbb{R}e (z)\bigr)^m \frac{(m - 1) !}{| \mathbb{R}e (2 z) |^m} .
  \label{eq:beta-exp-val}
\end{equation}
Due to the fact that usually $N \gg 1$, the minimal value of the right hand side of \eqref{eq:elleqerrmap} is attained for $m \gg 1$. 
So we expect that $\beta_{\ast}$ given in \eqref{eq:beta-exp} will be
a good choice that leads to good bounds for both $\| v - v_N \|^2$ and $| v -
v_N |_1^2$.

\begin{figure}[htbp]
  \centering
  \includegraphics[width=0.8\textwidth]{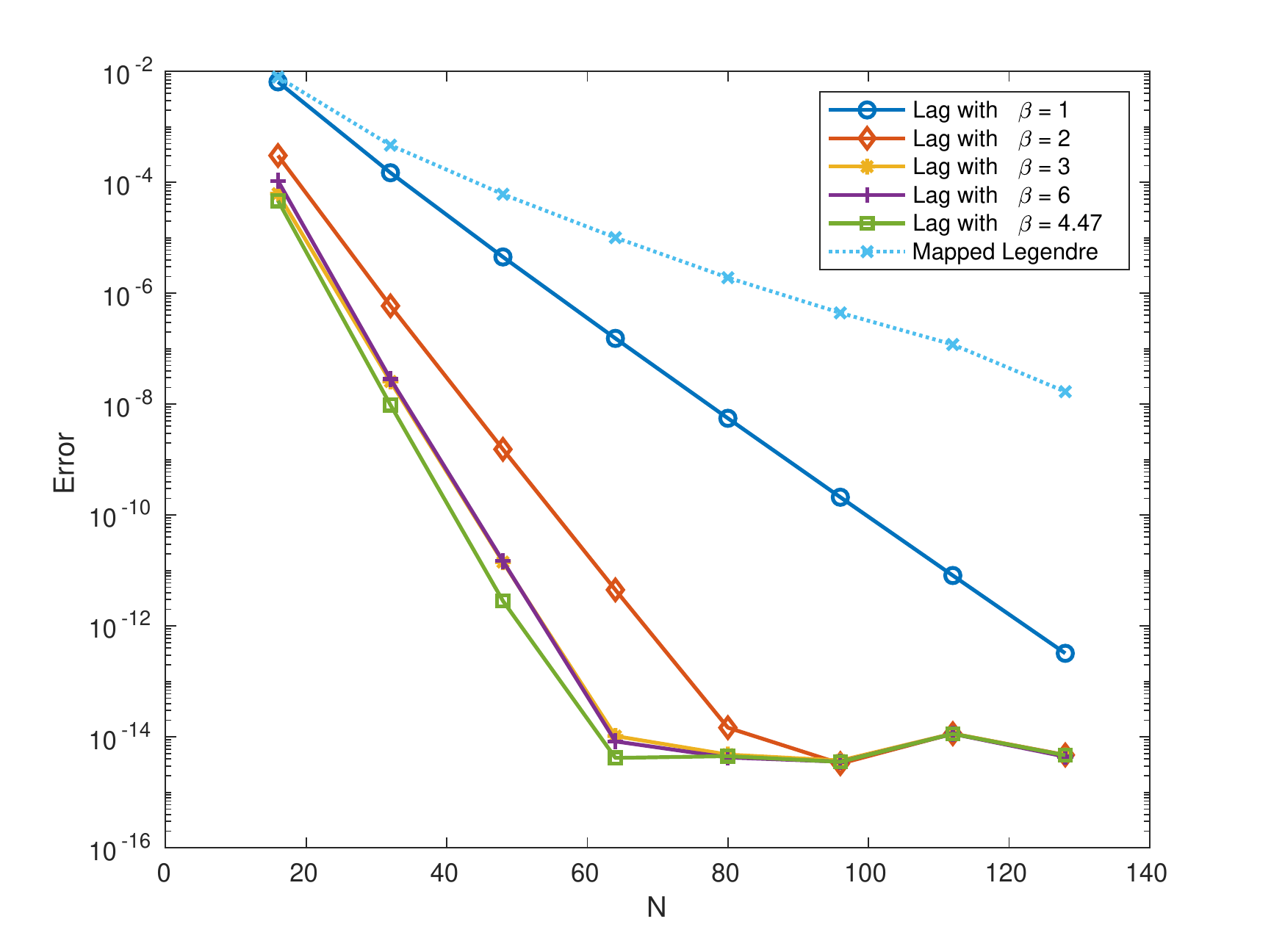}
  \caption{\label{fig:optscale}Numerical errors of the Laguerre method solving
  equation \eqref{eq:elleq} with $\gamma = 2$ and exact solution $u_1 (x) =
  \sin (2 x) e^{- x}$ using different values for scaling constant $\beta$, and
  compared with the mapped Legendre method.}
\end{figure}

In Figure \ref{fig:optscale}, we plot the numerical errors of the scaled
Laguerre method solving \eqref{eq:elleq} with an exact solution $u_1 (x) =
\sin (2 x) e^{- x}$ using different values of scaling factor $\beta$. To
compare its performance, the result of a mapped Legendre method (see
{\cite{shen_recent_2009}}) is also presented. The optimal value given by
\eqref{eq:beta-exp} is $\beta^{\ast} = 2 \sqrt{1 + 2^2} \cong 4.47$. We see
from the figure this choice of scaling indeed give best convergence speed. The
Laguerre method with $\beta$ close to $\beta^{\ast}$ produce results much
better than Laguerre method without scaling and the result of the mapped
Legendre method. We notice that the mapped Legendre method also allows a scaling
factor, but the dependence of convergence speed of mapped Legendre method on
the scaling factor is not as sensitive as the Laguerre method. By using
different scaling values, the mapped Legendre method can get convergence speed similar
to the Laguerre method with $\beta = 1$, but never close to the results of
Laguerre method with $\beta \sim \beta^{\ast}$.

For functions having other decaying properties, it is not easy to obtain an
explicit formula for the optimal value of the scaling factor $\beta$. We leave the problem of finding optimal scaling of the Laguerre method for general functions to a
future study.

\subsection{Using more than one thousand Laguerre bases}

In previous numerical studies using Laguerre methods for unbounded domains,
only small numbers of Laguerre bases are used due to the fact that the
standard methods to generate Laguerre polynomials and Laguerre functions
encounter serious round-off error and overflow/underflow issues, which are explained in last section. Now, we apply the
improved Laguerre algorithms to solve equation \eqref{eq:elleq} with
algebraic decay solutions $u_2$ and $u_3$. It is known that the convergence
rate is also algebraic~\cite{shen_stable_2000,shen_recent_2009}, thus, one may need a lot of Laguerre bases to get good
numerical results when the rate is slow.

\begin{figure}[htbp]
  \centering
  \includegraphics[width=0.8\textwidth]{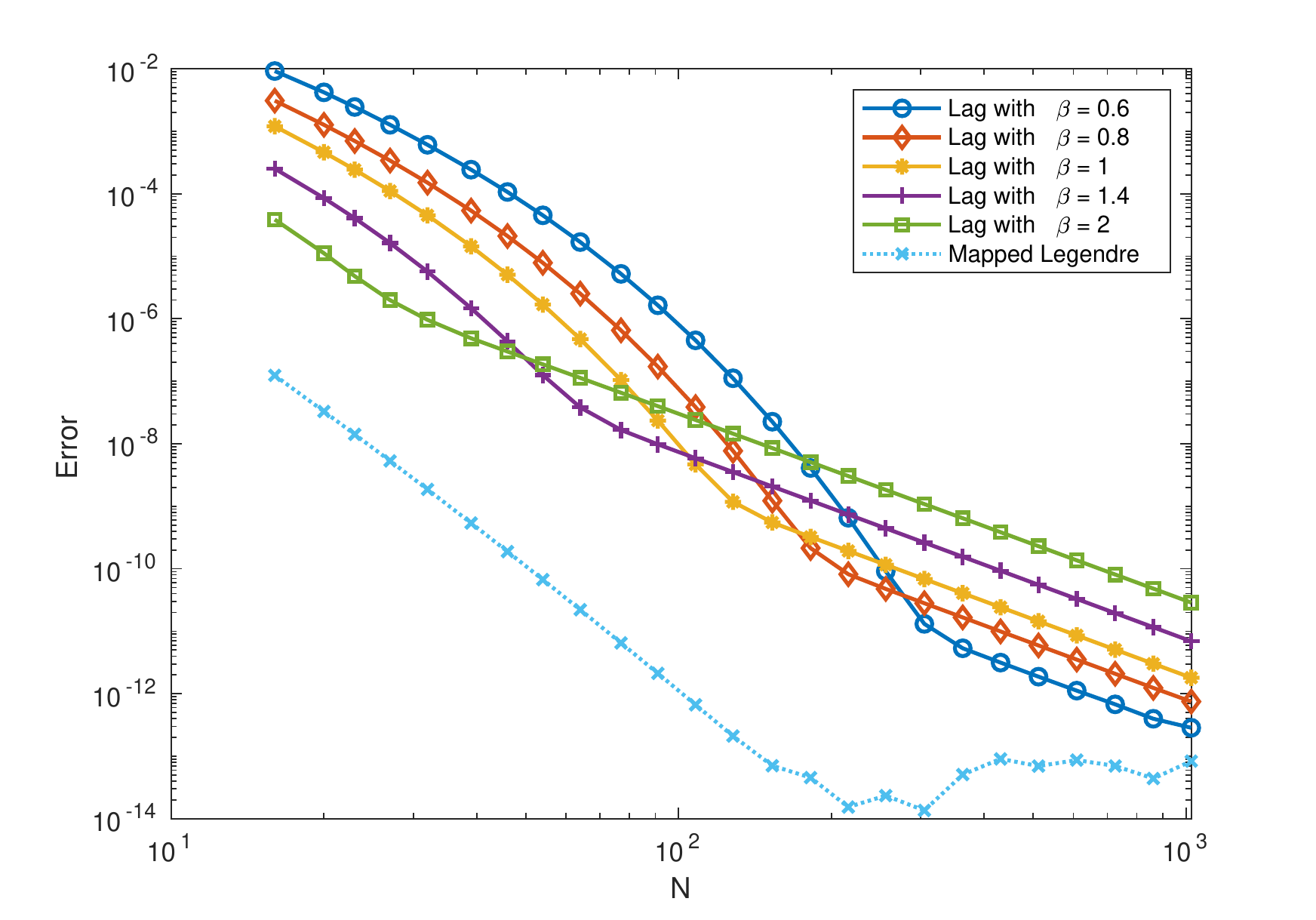}
  \caption{\label{fig:u2}Convergence results of the scaled Laguerre method and
  mapped Legendre method solving equation \eqref{eq:elleq} with $\gamma = 2$
  and exact solution $u_2 (x) = (1 + x)^{-5/2}$. The $L^2$ error is
  shown.}
\end{figure}

In Figure \ref{fig:u2}, we present the results of using scaled Laguerre
method and mapped Legendre method for \eqref{eq:elleq} with exact solution
$(1+x)^{-5/2}$. Note that, the tested case is same to the right plot of
Figure 8 in {\cite{shen_recent_2009}}, where only 128 bases are used. Here,
1024 Laguerre bases are used, which show more clear the convergence behavior
of the Laguerre method. In particular, Laguerre method with $\beta = 0.6$ and
1024 bases leads to $L^2$ error smaller than $3\times 10^{-13}$. From Figure
\ref{fig:u2}, we see if a smaller value of $\beta$ is used, one needs
use more degree of freedoms to pass through a larger pre-asymptotic range to
get the algebraic convergence rate. We also observe that larger scaling
factor $\beta$ should be used with less Laguerre bases, and smaller $\beta$
should be used with large numbers of Laguerre bases. This is
different to the exponential decay case. It can be observed  that the
converge of Laguerre method for this test case is slower than the mapped
Legendre method, this is consistent to the theoretical result~{\cite{shen_recent_2009}}.

\begin{figure}[htbp]
  \centering
  \includegraphics[width=0.8\textwidth]{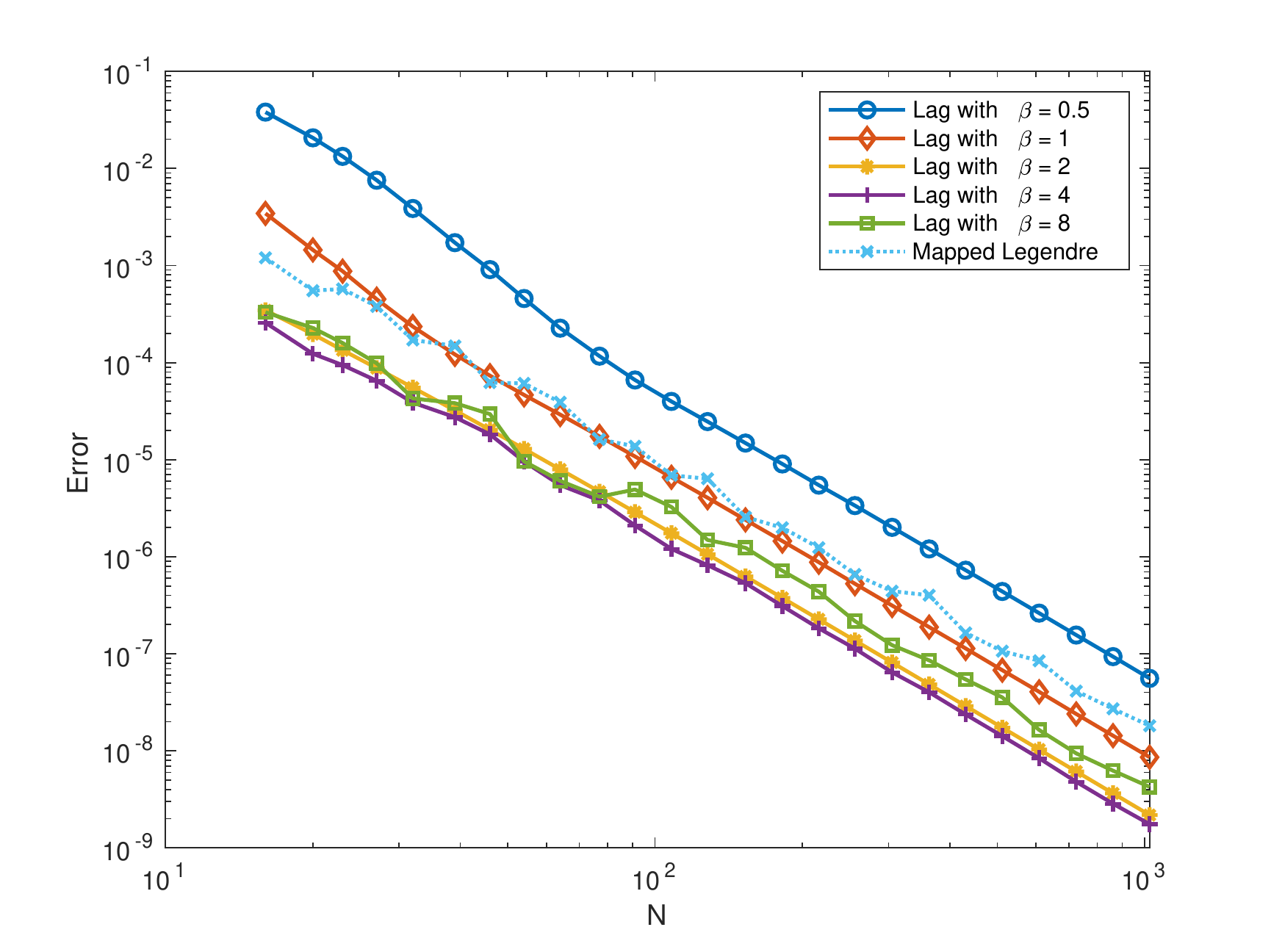}
  \caption{\label{fig:u3}Convergence results of the scaled Laguerre method and
  mapped Legendre method solving equation \eqref{eq:elleq} with $\gamma = 2$
  and exact solution $u_3 (x) = \sin(2x) (1+x)^{-7/2}$. The $L^2$
  error is shown.}
\end{figure}

In Figure \ref{fig:u3}, we show numerical results for algebraic decay
solutions with oscillations. Here, we take $u_3 (x) = \sin(2x)  (1 + x)^{-7/2}$. It is observed that the convergence speed of the Laguerre method with best scaling factor is
better than the mapped Legendre method. We also observe that the best
scaling factor is not sensitive to the number of bases used, which means that we can tune the scaling factor on coarse grids and use the result on refined grids. 

\section{Conclusion}

In this paper, we proposed improved algorithms to generate Laguerre polynomials,
Laguerre functions as well as Laguerre-Gauss quadrature schemes with better
stability and smaller round-off errors, which are essential to Laguerre
methods for solving partial differential equations in unbounded domains. These
algorithms can stably generate Laguerre functions of degree higher than one
thousand using floating-point number systems with much smaller numerical errors than existing methods. With the
help of these basic algorithms, one may use thousands of bases in the Laguerre
spectral methods to get better numerical results. We have demonstrated this by
considering the Laguerre spectral method solving an elliptic equation defined on
$\mathbb{R}_+$. We also presented some analysis and numerical investigations on the optimal scaling factor of the Laguerre spectral methods. We found that in the two special cases that 
the Laguerre method converges faster than the mapped Legendre method, the optimal 
scaling factor is independent of the numerical resolutions, which suggests that one can tune the scaling factor using coarse grids and use the results on refined grids.

\section*{Acknowledgments}
  We would like to thank Prof. Jie Shen,  Li-Lian Wang and Hui-Yuan Li for helpful discussions. This work was supported by the National Natural Science Foundation of China under Grant No. 12171467 and 11771439.
\bibliographystyle{siamplain}
\bibliography{Laguerre}

\end{document}